\documentclass[10pt,a4paper]{article}
\usepackage{latexsym,amssymb,amsmath,mathrsfs}
\usepackage{sectsty}
\usepackage{color}
\usepackage{bm}
\usepackage[utf8]{inputenc}
\usepackage[T1]{fontenc}
\usepackage[anchorcolor=blue,%
 bookmarks=true,%
 bookmarksnumbered=true,%
]{hyperref}
\usepackage{cite}
\usepackage{graphicx}
\allsectionsfont{\normalsize\sc\center}

\newtheorem{thm}{Theorem}[section]
\newtheorem{prop}[thm]{Proposition}
\newtheorem{cor}[thm]{Corollary}
\newtheorem{lem}[thm]{Lemma}
\newtheorem{defn}[thm]{Definition}
\newtheorem{rem}[thm]{Remark}
\newtheorem{ex}[thm]{Example}
\newtheorem{assumption}{(A)}

\newenvironment{pf}{\par\begin{trivlist}%
\item[]{\bf Proof.}\ }{\hfill $\square$ \end{trivlist}\par}
\newenvironment{apf}[1]{\par\begin{trivlist}%
\item[]{\bf Proof of #1.}\ }{\hfill $\square$ \end{trivlist}\par}

\makeatletter \@addtoreset{equation}{section} \makeatother

\newcommand{\R}{\mathbb{R}}

\newcommand{\N}{\mathbb{N}}

\DeclareMathOperator{\Ric}{Ric}

\renewcommand{\d}{\mathrm{d}}
\newcommand{\m}{\mathfrak{m} }

\newcommand{\PP}{\bf P}
\newcommand{\EE}{\bf E}

\newcommand{\E}{\mathscr{E}}
\newcommand{\F}{\mathscr{F}}

\newcommand{\1}{{\bf 1}}
\newcommand{\eps}{{\varepsilon}}

\newcommand{\<}{\langle}
\renewcommand{\>}{\rangle}

\title{\large\bf $L^p$-Kato class measures for symmetric Markov processes under heat kernel estimates  }
\author{Kazuhiro Kuwae\thanks{Department of Applied Mathematics, Fukuoka University,
Fukuoka 814-0180, Japan ({\sf kuwae@}{\sf fukuoka-u.ac.jp}). Supported in part by JSPS Grant-in-Aid for Scientific Research (KAKENHI) 17H02846 and by fund (No.:185001) from the Central Research Institute of Fukuoka University.}
\ \ and\ \ 
Takahiro Mori\thanks{Research Institute for Mathematical Sciences, Kyoto University, Kyoto, 606--8502, 
Japan ({\sf tmori@kurims.kyoto-u.ac.jp}). 
Supported by JSPS Grant-in-Aid for Scientific Research (KAKENHI) 18J21141. 
}
}
\date{}
\begin{document}
\maketitle
% \raggedbottom 

\begin{abstract}
In this paper, we establish the coincidence of two classes of
$L^p$-Kato class measures
in the framework of symmetric Markov processes admitting
upper and lower estimates of heat kernel
under mild conditions.
One class of $L^p$-Kato class measures is defined by
the $p$-th power of positive order resolvent kernel,
another is defined in terms of the $p$-th power of
Green kernel depending on some exponents
related to the heat kernel estimates.
We also prove that $q$-th integrable functions on balls
with radius $1$ having uniformity of its norm
with respect to centers are of $L^p$-Kato class
if $q$ is greater than a constant
related to $p$ and the constants
appeared in the upper and lower estimates of
 the heat kernel.
These are complete extensions of some results
by Aizenman-Simon
and the recent results by the second named author
in the framework of Brownian motions on Euclidean space.
We further give necessary and sufficient conditions
for a Radon measure with Ahlfors regularity
to belong to $L^p$-Kato class.
Our results can be applicable to many examples,
for instance, symmetric (relativistic) stable processes,
jump processes on $d$-sets,
Brownian motions on Riemannian manifolds,
diffusions on fractals and so on.
\end{abstract}

{\it Keywords}: Dirichlet form, Markov process, $L^p$-Kato class measure, $L^p$-Dynkin class measure, 
heat kernel, semigroup kernel, 
resolvent kernel, Green kernel, 
Stollman-Voigt inequality, 
Brownian motion, 
symmetric $\alpha$-stable process, 
relativistic $\alpha$-stable process, $d$-sets, Riemannian manifolds, 
Li-Yau's estimate, nested fractals,  Sierpi\'nski carpet

{\it Mathematics Subject Classification (2020)}: Primary 	60J25, 60J45, 60J46; Secondary 31C25, 35K08, 31E05.

\section{Introduction}\label{sec:introduction}
Let $p\in [1,+\infty[$.
A Borel measure  $\mu$ on $\R^d$ is said to be \emph{of $L^p$-Kato class} (\emph{of $p$-Kato class} in short) $K_d^{\,p}$ if
\begin{align*}
 \lim_{r\to0}
 \sup_{x\in\R^d}
 \int_{|x-y|<r}\frac{\,\mu(\d y)\,}{|x-y|^{(d-2)p}}
=
 0
\quad&
 \text{ for }\quad d\geq3
,
\\
 \lim_{r\to0}
 \sup_{x\in\R^d}
 \int_{|x-y|<r}(\log|x-y|^{-1})^p\mu(\d y)
=
 0
\quad&
 \text{ for }\quad d=2
,
\\
 \sup_{x\in\R^d}
 \int_{|x-y|\leq1}\mu(\d y)
<
 +\infty
\quad&
 \text{ for }\quad d=1
.
\end{align*}
We write $K_d$ instead of $K_d^1$ for $p=1$.
The notion of ($1$-)Kato class  $K_d$ was introduced by T.~Kato \cite{Kato:singular,Kato:vectorpoten}
in order to solve the essential self-adjointness of
the Schr\"odinger operator $-\Delta+V$ on $C_0^{\infty}(\R^d)$ (see the survey paper \cite{Simon:20} by Simon).
Let
$
 {\bf X}^{\rm w}
=
 (\Omega, B_t, {\PP}_x)_{x\in\R^d}
$
be a $d$-dimensional Brownian motion on $\R^d$.
The following theorem was proved
by Aizenman-Simon \cite{AizSim:Brown} under $p=1$ and noted
by the second named author in \cite{TM:pKato} for general $p\in]1,+\infty[$ with $d-p(d-2)>0$:

\begin{thm}[{\cite[Theorem~1.3(ii)]{AizSim:Brown}, \cite[Example~2.4]{TM:pKato}}]\label{thm:Kato}
Let $p\in[1,+\infty[$ with $d-p(d-2)>0$.
Then $\mu\in K_d^{\,p}$ if and only if
\begin{align*}
\sup\limits_{x\in\R^d}
 \int_{\R^d}
  \left(\int_0^tp_s(x,y)\d s\right)^p
\mu(\d y)
\longrightarrow0
\quad\text{ as }\quad t\to0,\label{eqn:Katosemigroup}
\end{align*}
where
$
 p_t(x,y)
:=
 \frac1{(2\pi t)^{d/2}}
 \exp\bigl(-\frac{\,|x-y|^2\,}{2t}\bigr)
$
is the heat kernel of ${\bf X}^{\rm w}$.
\end{thm}

\noindent
In \cite{KwT:Katounderheat},
we extend Theorem~\ref{thm:Kato} under $p=1$
in a probabilistic way.
Our first main theorem (Theorem \ref{Th2} below) is 
an extension of Theorem \ref{thm:Kato}
for general $p\in [1, +\infty[$
under the framework of \cite{KwT:Katounderheat}. 

\bigskip

The following theorem is proved
by Vogt~\cite[Proposition~2.2]{Vogt:Poincare2009}:

\begin{thm}[{Vogt~\cite[Proposition~2.2]{Vogt:Poincare2009}}]\label{thm:Vogt}
Let $\gamma\in]0,2]$. 
Suppose that a Borel measure $\mu$ on $\R^d$ satisfies that there exists $C>0$ such that 
$\mu(B_r(x))\leq Cr^{d-\gamma}$ for all $x\in \R^d$ and $r\in]0,+\infty[$. Then $\mu\in K_d^{\,p}$ with $p=1$ for $\gamma<2$ {\rm(}$\gamma\leq1$ if $d=1${\rm)}.  
\end{thm}

\noindent
Our second main theorem (Theorem \ref{Th5} below) is
a complete extension of Theorem~\ref{thm:Vogt}
for general $p\in[1,+\infty[$.

\bigskip

\par The following theorems are also shown by
Aizenman-Simon \cite{AizSim:Brown} under $p=1$ and noted in
\cite[Example~2.4]{TM:pKato} for general $p\in]1,+\infty[$
with $d-p(d-2)>0$:

\begin{thm}[{\cite[Theorem~1.4(iii)]{AizSim:Brown},\cite[Example~2.4]{TM:pKato}}]\label{thm:AizSim}
Let $p\in[1,+\infty[$ with $d-p(d-2)>0$.
Then $f\in L^q_{\text{\rm\tiny unif }}(\R^d)$ 
implies $|f|\d\m\in K_d^{\,p}$
if $q>d/(d-p(d-2))$ with $d\geq2$, or $q\geq1$ with $d=1$.
Here $\m$ denotes the Lebesgue measure on $\R^d$ and
$f\in L^q_{\text{\rm\tiny unif }}(\R^d)$ means
\begin{equation*}
\sup_{x\in\R^d}\int_{|x-y|\leq1}|f(y)|^q\d y<+\infty.
\end{equation*}
\end{thm}

\noindent
Let $M_{\alpha,q}$ $(\alpha>0,q>1)$
be the family of measurable functions $f$ on $\R^d$ satisfying
\begin{align*}
 \sup_{x\in \R^d}
 \int_{|x-y|\leq1}
  \frac{\;|f(y)|^q\,}{|x-y|^{d-\alpha}}
 \d y
<
 +\infty.
\end{align*}
$M_{\alpha,q}$ is called the \emph{Schechter class}.

\begin{thm}[{\cite[Propositions~4.1 and 4.2]{AizSim:Brown},\cite[Example~2.4]{TM:pKato}}]\label{thm:AizSimSchechter}
Let $p\in[1,+\infty[$ with $d-p(d-2)>0$.
Assume $d\geq3$.
Then 
$f\in M_{\alpha,q}$ implies $|f|\d\m\in K_d^{\,p}$
if $q>\alpha/(d-p(d-2))$.
\end{thm}

\begin{rem}
{\rm As noted in \cite[Example~2.4]{TM:pKato},
there are typos in \cite[Propositions~4.1 and 4.2]{AizSim:Brown}.
}
\end{rem}

\noindent
In \cite{KwT:Katofunc}, 
the first named author and Takahashi partially extend 
Theorem~\ref{thm:AizSim} for $p=1$ by
replacing $L^q_{\text{\rm\tiny unif }}(\R^d)$ (resp.~$K_d$)
with
$L^q(\R^d)$ (resp.~$S_K$) under
Nash-type estimate of semigroup kernel of Markov processes.
In \cite{KwT:Katounderheat}, 
they finally extend Theorem~\ref{thm:AizSim}
for $p=1$
in the framework of symmetric Markov processes
satisfying conditions
{\bf (A)}\ref{asmp:SC},
{\bf (A)}\ref{asmp:Bishop} and
{\bf (A)}\ref{asmp:Phi} below.
However, 
they did not write down the extension of
Theorem~\ref{thm:AizSimSchechter} under $p=1$
in \cite{KwT:Katounderheat}
for the limit of the length of \cite{KwT:Katounderheat}.

\bigskip

The purpose of this paper is to show that
the assertions in
Theorems~\ref{thm:Kato}, \ref{thm:Vogt}, \ref{thm:AizSim} and
\ref{thm:AizSimSchechter} can be extended for general $p\in[1,+\infty[$ in the framework of general symmetric
Markov processes admitting semigroup kernel with upper and lower estimates under some conditions.
These are presented as Theorems~\ref{Th2}, \ref{Th5}, \ref{Th3} and
\ref{Th4} in this paper.
Not only these extensions, we provide some useful
criteria for measures of $p$-Kato class in
Theorems~\ref{Th1}, \ref{Th6} and
Corollaries~\ref{cor:Th5}, \ref{cor:Th3}.
%In particular,
%any Hausdorff measure with its Hausdorff dimension 
%{\color{blue}{
%$\eta\in]\,(d-\alpha)_+,\,d\,]$ (resp.~$\eta\in]\,(d-2)_+,\,d\,]$)}}
%is of $p$-Kato class in the framework of
%rotationally symmetric $\alpha$-stable processes
%(resp.~Brownian motions on Riemmanian manifolds),
%provided {\color{blue}{$p\in[\,1,\,\eta/(d-\alpha)_+[$
%(resp.~$\eta\in[\,1,\,\eta/(d-2)_+[$)}}.
Our results are applicable to many Markov processes, for example, symmetric $\alpha$-stable
processes, relativistic $\alpha$-stable processes, jump type processes on $d$-sets,
Brownian motions on Riemannian manifolds with Ricci curvature lower bound and positivity of
injectivity radius,
diffusions on fractals and so on.

The constitution of this paper is as follows.
In Section~\ref{sec:asmp},
we prepare our framework and expose our assumptions. 
In Section~\ref{sec:maintheorem},
we state our main theorems
(Theorems~\ref{Th2}, \ref{Th5} and Corollary~\ref{cor:Th5}).
In Section~\ref{sec:criteria},
we state 
Theorems~\ref{Th1}, \ref{Th3}, \ref{Th4}, \ref{Th6}
and Corollary~\ref{cor:Th3}, which are useful criteria for 
$L^p$-Kato or $L^p$-Dynkin classes. 
Theorem~\ref{Th2} (resp.~{Theorems~\ref{Th5}, \ref{Th3}, \ref{Th4}) extends
Theorem~\ref{thm:Kato} (resp.~Theorems~\ref{thm:Vogt},
\ref{thm:AizSim}, \ref{thm:AizSimSchechter}).
In Section~\ref{sec:proof},
we give the proofs of Theorems~\ref{Th2}, \ref{Th5}, \ref{Th1}, \ref{Th3}, \ref{Th4},  \ref{Th6} and Corollaries~\ref{cor:Th5}, \ref{cor:Th3}.
In the last section, we expose concrete examples.

Finally, we announce the content of the paper \cite{KwM:pGreentightKato}
on the $L^p$-Green-tight measures of Kato class.
In \cite{KwM:pGreentightKato}, for transient symmetric Markov processes, we investigate the family
$S_{C\!K_{\infty}}^{\,p}$\!\! of $L^p$-Green-tight measures of $L^p$-Kato class in the sense of Z.-Q.~Chen, which is defined to be a subclass of
$S_K^{\,p}$ and
under the global heat kernel estimate, we 
prove the coincidence with the family 
$K_{\nu,\beta}^{p, \infty}$ of $L^p$-Green-tight measures of $L^p$-Kato class in the sense of Zhao which is defined to be a subclass of $K_{\nu,\beta}^{\,p}$. This is a natural extension of our
Theorem~\ref{Th2}.

\section{Preliminary}\label{sec:asmp}
For real numbers $a,b\in\R$, we set $a\lor b:=\max\{a,b\}$ and $a\land b:=\min\{a,b\}$.
Let $(E,d)$ be a locally compact separable metric space and $\m$ a positive Radon measure with full support.
Let $E_{\partial}:=E\cup\{\partial\}$ be 
the one-point compactification of $E$. 
For each $x\in E$ and $r>0$, denote by $B_r(x):=\{y\in E\mid d(x,y)<r\}$
the open ball with center $x$ and radius $r$.
We consider and fix a symmetric regular
Dirichlet form $(\E,\F)$ on $L^2(E;\m)$. Then there exists a Hunt process
 ${\bf X}=(\Omega,X_t,\zeta,{\PP}_x)$ such that for each Borel $u\in L^2(E;\m)$,
$T_tu(x)={\EE}_x[u(X_t)]$ $\m$-a.e.~$x\in E$ for all $t>0$, where  $(T_t)_{t>0}$
is the
semigroup associated with $(\E,\F)$. Here $\zeta:=\inf\{t\geq0\mid X_t=\partial\}$
denotes the life time of {\bf X}. For a Borel set $B$,
we denote $\sigma_B:=\inf\{t>0\mid X_t\in B\}$ (resp.~ $\tau_B:=\inf\{t>0\mid X_t\notin B\}$)
the \emph{first hitting time to $B$} (resp.~\emph{first exit time from $B$}).
Further, we assume that there exists a jointly measurable function $p_t(x,y)$ defined for all $(t,x,y)\in
]0,+\infty[\times E\times E$ such that ${\EE}_x[u(X_t)]=P_tu(x):=\int_Ep_t(x,y)u(y)\m(\d y)$
for any $x\in E$,
bounded Borel  function $u$ and $t>0$.
$p_t(x,y)$ is said to be a \emph{semigroup kernel}, or  sometimes called a \emph{heat kernel} of {\bf X}
on the analogy of heat kernel of diffusions.
Then $P_t$ can be extended to contractive semigroups on $L^p(E;\m)$ for $p\geq1$.
The following are well-known:
\begin{enumerate}
\item[(1)]\label{item:kernel1}
$
 \displaystyle
 p_{t+s}(x,y)=\int_E p_{s}(x,z)p_{t}(z,y)\m(\d z)
$
\quad
for \ all \quad $x,y \in E$ \quad and \quad $t,s > 0$.
\item[(2)]\label{item:kernel2}
$P_{t}(x,\d y)=p_{t}(x,y)\m(\d y)$
\quad
for \ all \quad $ x \in E$\quad and\quad $ t > 0$.
\item[(3)]\label{item:kernel3}
$
 \displaystyle
 \int_E p_{t}(x,y)\m(\d y) \leq 1
$
\quad
for \ all \quad $x \in E$\quad and\quad $t > 0$.
\end{enumerate}

Throughout this paper, we fix $\nu,\beta\in]0,+\infty[$ and $t_0\in]0,+\infty]$ and prepare
the following assumptions.

\begin{assumption}[Life time condition]\label{asmp:SC}
{\bf X} has the following property that
\begin{equation*}
\lim_{t\to0}\sup_{x\in E}{\PP}_x(\zeta\leq t)=:\gamma\in[0,1[.
\end{equation*}
In particular, if {\bf X}
is stochastically complete, that is, {\bf X} is conservative, then this condition is satisfied
with $\gamma=0$.
\end{assumption}

We fix an increasing positive function $V$ on $]0,+\infty[$.
\begin{assumption}[Bishop type inequality]\label{asmp:Bishop}
Suppose $r\mapsto V(r)/r^{\nu}$ is increasing or bounded, and
$\sup\limits_{x\in E}\m(B_r(x))\leq V(r)$ for all $r>0$.
\end{assumption}

\begin{assumption}[Upper and lower estimates of heat kernel]\label{asmp:Phi}
Let $\Phi_i$  $(i=1,2)$  be positive decreasing  functions
defined on $[0,+\infty[$ which may depend on $t_0$ if $t_0<+\infty$
and assume that $\Phi_2$ satisfies the following condition $H(\Phi_2)$:
\begin{equation*}
\int_1^{\infty}\frac{\,(V(t)\lor t^{\nu})\Phi_2(t)\,}{t}\d t<+\infty
\end{equation*}
and $(\Phi {\EE}_{\nu,\beta})$: for any
$x,y\in E$, $t\in]0,t_0[$
\begin{align*}
\frac{\,1\,}{t^{\nu/\beta}}\Phi_1\left(\frac{\,d(x,y)\,}{t^{1/\beta}}\right)
\leq p_t(x,y)\leq
\frac{\,1\,}{t^{\nu/\beta}}\Phi_2\left(\frac{\,d(x,y)\,}{t^{1/\beta}}\right).
\end{align*}
\end{assumption}
Note that the assumption {\bf (A)}\ref{asmp:Phi}
is essentially introduced 
as
the hypothesis of \cite{Gri:heatMMS}.

\bigskip

We next introduce
the classes of measures dealt with in this paper.
Throughout this paper, we consider a constant $p\in[1,+\infty[$.
\begin{defn}[{\boldmath$L^p$-Kato class \boldmath$S_K^{\,p}$, \boldmath$L^p$-Dynkin class \boldmath$S_D^{\,p}$}]
{\rm For a positive Radon measure $\mu$ on $E$, $\mu$ is said to be \emph{of $L^p$-Kato {\rm(}$p$-Kato in short{\rm)} class
relative to $p_t(x,y)$} (write $\mu\in S_K^{\,p}$) if
\begin{align}
\lim_{t\to0}\sup_{x\in E}\int_E
\left(\int_0^t p_s(x,y)\d s\right)^p
\mu(\d y)=0\label{SK1}
\end{align}
and $\mu$ is said to be
\emph{of $L^p$-Kato {\rm(}$p$-Kato in short{\rm)} class
relative to $p_t(x,y)$ of order $\delta\in]0,1]$} (write $\mu\in S_K^{\,p,\delta}$) if
\begin{align}
\sup_{x\in E}\left(\int_E\left(\int_0^tp_s(x,y)\d s\right)^p\mu(\d y)\right)^{\frac{\,1\,}{p}}=O(t^{\delta})\quad (t\to0).\label{SK1order}
\end{align}
For a positive Radon measure $\mu$ on $E$,
$\mu$ is said to be \emph{of $L^p$-Dynkin {\rm(}$p$-Dynkin in short{\rm)} class
relative to $p_t(x,y)$} (write $\mu\in S_D^{\,p}$) if
\begin{align}
 \sup_{x\in E}
 \int_E\left(\int_0^tp_s(x,y)\d s\right)^p\mu(\d y)
<
 +\infty
\quad \text{for \ some } t>0.\label{SD1}
\end{align}
Clearly, $S_K^{\,p,\delta}\subset S_K^{\,p}\subset S_D^{\,p}$.
When $p=1$, we write $S_D$ (resp.~$S_K$,
$S_K^{\,\delta}$
) instead of
$S_D^1$ (resp.~$S_K^1$, $S_K^{\,1,\delta}$
) for simplicity.
}
\end{defn}

For $\alpha>0$,
we denote the {\it $\alpha$-order resolvent kernel} by 
$
 r_{\alpha}(x,y):=\int_0^{\infty}e^{-\alpha t}p_t(x,y)\d t
$.
The following are proved by the second named author
in \cite[Proposition~2.6 and Corollary~2.7]{TM:pKato}
extending \cite[Lemma~3.1]{KwT:Katofunc}.

\begin{lem}[{\cite[Propositions~2.6, 2.10 and Corollary~2.7]{TM:pKato}
}]\label{lem:equiKatoDynkin}
Let $\delta\in]0,1]$.
$\mu\in S_K^{\,p}$ {\rm(}resp.~$\mu\in S_K^{\,p,\delta}${\rm)} is equivalent to
\begin{align}
\lim_{\alpha\to\infty}\sup_{x\in E}\int_E r_{\alpha}(x,y)^p\mu(\d y)=0\label{SK2}
\end{align}
\begin{align}
\text{\rm (resp.~}\quad \sup_{x\in E}\left(\int_E r_{\alpha}(x,y)^p\mu(\d y)\right)^{\frac{\,1\,}{p}}=O(\alpha^{-\delta})\quad (\alpha\to\infty)\quad\text{\rm )}\label{SK2order}
\end{align}
and $\mu\in S_D^{\,p}$ is equivalent to
\begin{align}
 \sup_{x\in E}
 \int_E r_{\alpha}(x,y)^p\mu(\d y)
<
 +\infty
\quad \text{for \ some }\alpha>0.\label{SD2}
\end{align}
\end{lem}

\begin{lem}[{\cite[Proposition~2.6]{TM:pKato}, see also \cite[Lemma~3.2]{KwT:Katofunc},\cite{AM:AF}}]\label{lem:Dynkinequi}
The following are equivalent to each other: 
\begin{enumerate}
\item\label{item:SD}
$\mu\in S_D^{\,p}$.
\item\label{item:PT}
$
 \displaystyle
 \sup\limits_{x\in E}
 \int_E\left(\int_0^tp_s(x,y)\d s\right)^p\mu(\d y)
<
 +\infty
\quad \text{for \ any }\quad t>0
$.
\item\label{item:RA}
$
 \displaystyle
 \sup\limits_{x\in E}
 \int_E r_{\alpha}(x,y)^p\mu(\d y)
<
 +\infty
\quad \text{for \ any }\quad \alpha>0
$.
\end{enumerate}
\end{lem}

\begin{defn}[Dynkin class \boldmath$D_{\nu,\beta}^{\,p}$]\label{df:DynkinGreen}
{\rm
Fix $\nu>0$ and $\beta>0$.
For a positive Borel measure $\mu$ on $E$, $\mu$ is said to be \emph{of $L^p$-Dynkin ($p$-Dynkin in short) class
relative to Green kernel} (write $\mu\in D_{\nu,\beta}^{\,p}$) if 
\begin{align*}
 \sup_{x\in E}
 \int_{d(x,y)<r}G(x,y)^p\mu(\d y)
&<
 +\infty
\quad
\text{ for \ some }\quad r>0
\quad
 \text{for }\nu\geq\beta,
\\
 \sup_{x\in E}
 \int_{d(x,y)\leq 1}\mu(\d y)
&<
 +\infty
\quad
 \text{ for }\quad \nu<\beta,
\end{align*}
where $G(x,y):=G(d(x,y))$ with 
\begin{align*}
G(r) := 
\left\{\begin{array}{ll}
                     r^{\beta - \nu} &  \nu> \beta, \quad r\in]0,+\infty[,\\
                     \log(r^{-1}) &  \nu=\beta, \quad r\in]0,1[.
                     \end{array}\right.
\end{align*}
When $p=1$, we write $D_{\nu,\beta}$ 
instead of 
$D_{\nu,\beta}^{\,1}$.  
}
\end{defn}

\begin{defn}[Kato class \boldmath$K_{\nu,\beta}^{\,p}$]\label{df:KatoGreen} 
{\rm
Fix $\nu>0$ and $\beta>0$. 
For a positive Borel measure $\mu$ on $E$, $\mu$ is said to be \emph{of $L^p$-Kato ($p$-Kato in short) class 
relative to Green kernel} (write $\mu\in K_{\nu,\beta}^{\,p}$) if 
\begin{align*}
 \lim_{r\to0}
 \sup_{x\in E}&
 \int_{d(x,y)<r}G(x,y)^p\mu(\d y)
=
 0
\quad
 \text{ for  }\quad \nu\geq\beta,
\\
 \sup_{x\in E}&
 \int_{d(x,y)\leq 1}\mu(\d y)
<
 +\infty
\quad
 \text{ for }\quad \nu<\beta,
\end{align*}
where $G(x,y)$ is the function appeared above.  
When $p=1$, we write $K_{\nu,\beta}$ 
instead of 
$K_{\nu,\beta}^{\,1}$. 
}
\end{defn}

\begin{lem}\label{lem:KatoRadon}
If $\mu\in D_{\nu,\beta}^{\,p}$, then $\sup_{x\in E}\mu(B_r(x))<+\infty$ for small $r\in]0,e^{-1}[$. 
In particular, every $\mu\in D_{\nu,\beta}^{\,p}
$ is a Radon measure. 
\end{lem}

\begin{pf} The assertion is clear from $\mu(B_r(x))\leq\frac{\,1\,}{G(r)^p}\int_{B_r(x)}G(x,y)^p\mu(\d y)$ for 
$r\in]0,e^{-1}[$ with $\nu\geq\beta$. The case for $\nu<\beta$ is trivial.   
\end{pf}

\begin{lem}\label{lem:monotonicityKato}
For $1\leq p_1\leq p_2$, we have 
$D_{\nu,\beta}^{\,p_2}\subset D_{\nu,\beta}^{\,p_1}$ and 
$K_{\nu,\beta}^{\,p_2}\subset K_{\nu,\beta}^{\,p_1}$.
In particular,
$D_{\nu,\beta}^{\,p}\subset D_{\nu,\beta}$ and 
$K_{\nu,\beta}^{\,p}\subset K_{\nu,\beta}$ hold.  
\end{lem}

\begin{pf}
When $\nu<\beta$, $K_{\nu,\beta}^{\,p}=D_{\nu,\beta}$ is independent of $p$. 
So we may assume $\nu\geq\beta$. 
Let $1\leq p_1\leq p_2$ and take $\mu\in D_{\nu,\beta}^{\,p_2}$. Then
\begin{align*}
\int_{B_r(x)}G(x,y)^{p_1}\mu(\d y)&\leq 
\left(\int_{B_r(x)}G(x,y)^{p_2}\mu(\d y) \right)^{\frac{\,p_1\,}{p_2}}
\left(\mu(B_r(x)) \right)^{1-\frac{\,p_1\,}{p_2}}\\
&\leq \left(\int_{B_r(x)}G(x,y)^{p_2}\mu(\d y) \right)^{\frac{\,p_1\,}{p_2}}
\left(\frac{\,1\,}{G(r)^{p_2}}\int_{B_r(x)}G(x,y)^{p_2}\mu(\d y) \right)^{1-\frac{\,p_1\,}{p_2}}\\
&= \frac{\,1\,}{G(r)^{p_2-p_1}}\int_{B_r(x)}G(x,y)^{p_2}\mu(\d y),\quad r\in]0,e^{-1}[,
\end{align*}
which 
implies $\nu\in D_{\nu,\beta}^{\,p_2}$. 
Since $\lim_{r\to0}1/G(r)=0$, we obtain the 
inclusion $K_{\nu,\beta}^{\,p_2}\subset K_{\nu,\beta}^{\,p_1}$. 
\end{pf}

\begin{defn}[Measures of finite energy integrals, \boldmath$S_0$, \boldmath$S_{00}$; cf.~\cite{FOT}]\label{df:S_0,S_00}
{\rm
A Radon measure $\mu$ on $E$ is said to be
\emph{of finite energy integral} with respect to $(\E,\F)$
(write $\mu\in S_0$) if there exists $C>0$ such that
\begin{equation*}
 \int_E|v|\d\mu
\leq
 C\sqrt{\E_1(v,v)},
\quad
 \text{for \ any }\quad v\in\F\cap C_0(E).
\end{equation*}
In that case, for every $\alpha >0$, there exists 
$U_{\alpha }\mu\in \F$ such that 
\begin{equation*}
 \E_{\alpha}(U_{\alpha }\mu,v)
=
 \int_E v(x)\mu(\d x),
\quad
 \text{for \ any }\quad v \in \F \cap C_0(E).
\end{equation*} 
Moreover we write $\mu\in S_{00}$ if
$\mu(E)<+\infty$ and $U_{\alpha}\mu\in\F\cap L^{\infty}(E;\m)$
for some/all $\alpha>0$. 
}
\end{defn}

\begin{defn}[Smooth measures, \boldmath$S$; cf.~\cite{FOT}]\label{df:smooth}
{\rm
A Borel measure $\mu$ on $E$ is said to be a
\emph{smooth measure} with respect to $(\E,\F)$ 
(write $\mu\in S$) if $\mu$ charges no exceptional set and 
there exists a generalized nest $\{F_n\}$ of closed sets
such that $\mu(F_n)<+\infty$ for each $n\in\N$. 
}
\end{defn}

\begin{defn}[Smooth measures in the strict sense, \boldmath$S_1$; cf.~\cite{FOT}]\label{df:smoothstrict}
{\rm 
A Borel measure $\mu$ on $E$ is said to be a \emph{smooth measure in the strict sense} with respect to $(\E,\F)$ (write $\mu\in S_1$) if there exists an increasing sequence $\{E_n\}$ of Borel sets such that 
$E = \bigcup_{n=1}^{\infty}E_n$, and for any $n\in \N$, 
$\1_{E_n}\mu \in S_{00}$ and ${\PP}_x(\lim_{n \to \infty }\sigma_{E\setminus E_n}\geq \zeta)=1$  
 for any $x\in E$. 
}
\end{defn}

\begin{rem}
{\rm 
It is shown in \cite[Proposition~2.5]{TM:pKato} that $S_K^{\,p}\subset S_D^{\,p}\subset S_1$. 
}
\end{rem}

\section{Main theorems}\label{sec:maintheorem}

Now we are ready to state the main theorems.
Our first main theorem 
is a complete extension of 
Theorem~\ref{thm:Kato} and 
\cite[Theorem~3.2]{KwT:Katounderheat}. 
This is the most important theorem in this paper. 

\begin{thm}\label{Th2}
Let $p\in[1,+\infty[$.  
Suppose that  {\bf (A)}\ref{asmp:SC}, {\bf (A)}\ref{asmp:Bishop} and {\bf (A)}\ref{asmp:Phi} hold. 
Then we have $S_K^{\,p}=K_{\nu,\beta}^{\,p}$ and 
$S_D^{\,p}=D_{\nu,\beta}^{\,p}$.  Moreover, 
$\mu\in K_{\nu,\beta}^{\,p}$ implies that 
\begin{align}
\sup_{x\in E}\mu(B_R(x))<+\infty\quad\text{ for \ all }\quad R>0.\label{eq:ballunifinite}
\end{align} 
For $\nu<\beta$, we have 
$S_D^{\,p}=K_{\nu,\beta}^{\,p}=S_D=K_{\nu,\beta}$ and $\mu\in K_{\nu,\beta}^{\,p}$ is equivalent to 
(\ref{eq:ballunifinite}). 
\end{thm}

Our second main theorem
gives a criterion for 
$L^p$-Kato and $L^p$-Dynkin class measures 
based on the decay rate of the measures of balls. 

\begin{thm}\label{Th5} 
Let $\mu$ be a Radon measure,
$p\in[1,+\infty[$ and $\eta \in ]0, \nu]$.
Suppose that
{\bf (A)}\ref{asmp:Bishop} and {\bf (A)}\ref{asmp:Phi} hold. 
\begin{enumerate}
\item\label{item:Th5:1}
If there exist constants $r_0, C_2>0$
such that 
$\mu(B_r(x))\leq C_2 r^{\eta}$
for any $x\in E$ and $r\in]0,r_0]$ and
$\eta-p(\nu-\beta)>0$ holds,
then $\mu\in K_{\nu,\beta}^{\,p}$. 
\item\label{item:Th5:2} 
If there exist $x_0\in E$ and constants $r_0, C_1>0$
such that 
$C_1 r^{\eta}\leq \mu(B_r(x_0))$
for any $r\in]0,r_0]$ and
$\mu\in D_{\nu,\beta}^{\,p}$
holds,
then $\eta-p(\nu-\beta)\geq0$.
\item\label{item:Th5:3} 
If there exist $x_0\in E$ and constants $r_0, C_1,C_2>0$
such that
$C_1 r^{\eta}\leq \mu(B_r(x_0))\leq C_2r^{\eta}$
for any $r\in]0,r_0]$ and
$\mu\in D_{\nu,\beta}^{\,p}$ holds,
then $\eta-p(\nu-\beta)> 0$.
\end{enumerate}
In particular,
if $\mu$ satisfies the Ahlfors regularity, i.e.,
$C_1r^{\eta}\leq\mu(B_r(x))\leq C_2r^{\eta}$
for all $x\in E$ and $r\in ]0,r_0]$ with some $r_0,C_1,C_2>0$,
then the following are equivalent:
\begin{enumerate}
\item[{\rm (1)}] $\mu\in K_{\nu,\beta}^{\,p}$.
\qquad
{\rm (2)} $\mu\in D_{\nu,\beta}^{\,p}$.
\qquad
{\rm (3)}  $\eta-p(\nu-\beta)>0$.
\end{enumerate}
\end{thm}

\begin{cor}\label{cor:Th5}
Let $p\in[1,+\infty[$.
Suppose that
{\bf (A)}\ref{asmp:SC},
{\bf (A)}\ref{asmp:Bishop} and
{\bf (A)}\ref{asmp:Phi} hold.  
Then the following are equivalent:
\begin{enumerate}
\item    $\m\in K_{\nu,\beta}^{\,p}=S_K^{\,p}$.
\qquad
{\rm(2)} $\m\in D_{\nu,\beta}^{\,p}=S_D^{\,p}$.
\qquad
{\rm(3)} $\nu-p(\nu-\beta)>0$.
\end{enumerate}
\end{cor}

\section{Criteria for $p$-Kato and $p$-Dynkin classes}\label{sec:criteria}

In this section,
we give other criteria for $p$-Kato and $p$-Dynkin classes. 

\begin{thm}\label{Th1} 
Let $p\in[1,+\infty[$. 
Suppose that {\bf (A)}\ref{asmp:Phi} 
and $\nu\geq\beta$ hold. Then the following are equivalent:
\begin{enumerate}
\item\label{item:Green} $\mu\in K_{\nu,\beta}^{\,p}$.
\item\label{item:resolventany} 
For any $\alpha>0$, \quad $\displaystyle{\lim_{r\to0}\sup_{x\in E}\int_{B_r(x)}r_{\alpha}(x,y)^p\mu(\d y)=0}$.
\item\label{item:resolventsome} 
For some $\alpha>0$, \quad $\displaystyle{\lim_{r\to0}\sup_{x\in E}\int_{B_r(x)}r_{\alpha}(x,y)^p\mu(\d y)=0}$. \item\label{item:heatany} 
For any $t>0$, \quad
$\displaystyle{\lim_{r\to0}\sup_{x\in E}\int_{B_r(x)}\left(\int_0^tp_s(x,y)\d s\right)^p\mu(\d y)=0}$. \item\label{item:heatsome} 
For some $t>0$, \quad
$\displaystyle{\lim_{r\to0}\sup_{x\in E}\int_{B_r(x)}\left(\int_0^tp_s(x,y)\d s\right)^p\mu(\d y)=0}$. \end{enumerate}
Moreover, the following are equivalent: 
\begin{enumerate}
\item[{\rm $(1')$}]\label{item:Green*}
$\mu\in D_{\nu,\beta}^{\,p}$.
\item[{\rm $(2')$}]\label{item:resolventany*} 
For any $\alpha>0$, \quad
$
 \displaystyle
 \sup_{x\in E}
 \int_{B_r(x)}r_{\alpha}(x,y)^p\mu(\d y)
<
 +\infty
$
\quad for some \quad $r>0$. 
\item[{\rm $(3')$}]\label{item:resolventsome*} 
For some $\alpha>0$, \quad
$
 \displaystyle
 \sup_{x\in E}
 \int_{B_r(x)}r_{\alpha}(x,y)^p\mu(\d y)
<
 +\infty
$ 
\quad for some \quad $r>0$. 
\item[{\rm $(4')$}]\label{item:heatany*} 
For any $t>0$, \quad
$
 \displaystyle
 \sup_{x\in E}
 \int_{B_r(x)}\left(\int_0^tp_s(x,y)\d s\right)^p\mu(\d y)
<
 +\infty
$ 
\quad for some \quad $r>0$.  
\item[{\rm $(5')$}]\label{item:heatsome*} 
For some $t>0$, \quad
$
 \displaystyle
 \sup_{x\in E}
 \int_{B_r(x)}\left(\int_0^tp_s(x,y)\d s\right)^p\mu(\d y)
<
 +\infty
$ 
\quad for some \quad $r>0$. 
\end{enumerate}
\end{thm}

\begin{rem}\label{rem:subproc} 
{\rm Theorem~\ref{Th1} is a complete extension of 
\cite[Theorem~3.1]{KwT:Katounderheat}.  
The equivalence among \ref{item:Green}--\ref{item:heatsome} in 
Theorem~\ref{Th1} does not hold for $\nu<\beta$ in general. 
In fact, for $1$-dimensional Brownian motion
${\bf X}^{\rm w}$,
we see that 
$\mu=\delta_0\in K_1=K_1^{\,p}$ does not satisfy the conditions \ref{item:resolventany}, \ref{item:resolventsome} in Theorem~\ref{Th1} because of  
$r_{\alpha}(x,y)=e^{-\sqrt{2\alpha}|x-y|}/\sqrt{2\alpha}$.
}
\end{rem}

The following theorem is a complete extension of 
Theorem~\ref{thm:AizSim} and
\cite[Theorem~3.3]{KwT:Katounderheat}.
Though this is a special case of 
Theorem~\ref{Th4} below, we expose it in connection with the previous results. 

\begin{thm}\label{Th3}
Let $p\in[1,+\infty[$ with $\nu-p(\nu-\beta)>0$. 
Suppose that {\bf (A)}\ref{asmp:Bishop} holds.
Then for any $f \in  L^q_{\text{\rm\tiny unif}}(E;\m)$, 
we have $|f|\d\m\in K_{\nu,\beta}^{\,p}$ if  $q > \nu/(\nu-p(\nu-\beta))$ with $\nu \geq \beta$, or if $q \geq 1$ 
with $\nu < \beta$.  
Here $f \in  L^q_{\text{\rm \tiny unif}}(E;\m)$ 
means 
\begin{equation*}
\sup\limits_{x\in E}\int_{d(x,y)\leq1}|f(y)|^q\m(\d y)<+\infty.
\end{equation*}
\end{thm} 

The following corollary is an easy consequence of 
Theorem~\ref{Th3}.

\begin{cor}\label{cor:Th3}
Let $p\in[1,+\infty[$ with $\nu-p(\nu-\beta)>0$. 
Suppose that {\bf (A)}\ref{asmp:Bishop} holds. 
For any fixed point $o\in E$, 
$d(\cdot,o)^{-\gamma}\m\in K_{\nu,\beta}^{\,p}$ if 
$\gamma\in[0,\nu-p(\nu-\beta)[$ with $\nu\geq\beta$, and $\gamma\in[0,1[$ with $\nu<\beta$. In particular, 
$\m\in K_{\nu,\beta}^{\,p}$ always holds. 
\end{cor}

Next theorem is a generalization of
\cite[Propositions~4.1 and 4.2]{AizSim:Brown},
which does not treat the case $q=1$.
The assertion of Theorem~\ref{Th3}
can be recovered from Theorem~\ref{Th4}
by setting $\nu=\alpha$.

\begin{thm}\label{Th4}
Let $p\in[1,+\infty[$ with $\nu-p(\nu-\beta)>0$. 
Suppose that {\bf (A)}\ref{asmp:Bishop} holds,
 and assume
$\nu \geq \alpha$ and $q\geq 1$, or
$\nu<\alpha$ and $q>\alpha/\nu$.
Then for any $f \in  M_{\alpha,q}$, 
we have $|f|\d\m\in K_{\nu,\beta}^{\,p}$ if  $q > \alpha/(\nu-p(\nu-\beta))$,
where $f \in  M_{\alpha,q}$ 
means
\begin{equation*}
\sup\limits_{x\in E}\int_{d(x,y)\leq1}\frac{\,|f(y)|^q\,}{d(x,y)^{\nu-\alpha}}\m(\d y)<+\infty.
\end{equation*}
\end{thm}

Finally, we give the following theorem without assuming the lower estimate of the heat kernel. 

\begin{thm}\label{Th6}
Let $p\in[1,+\infty[$.
Assume that  
there exist constants $C,r_0>0$ and $\eta \in]0,\nu]$
satisfying $\eta-p(\nu-\beta)>0$
such that 
$\mu(B_r(x))\leq Cr^{\eta}$
for any $x\in E$ and $r\in]0,r_0[$. 
Suppose that
{\bf (A)}\ref{asmp:Phi} holds without assuming 
the lower estimate of $(\Phi {\EE}_{\nu,\beta})$. 
More precisely,
there exists a positive decreasing function $\Phi_2$ on 
$[0,+\infty[$ which may depend on $t_0$ if $t_0<+\infty$ 
and 
$\Phi_2$ satisfies 
\begin{align*}
 \int_1^{\infty}t^{\nu-1}\Phi_2(t)\d t<+\infty
\end{align*}
and for any $x,y\in E$ and $t\in]0,t_0[$ 
\begin{align*}
 p_t(x,y)
\leq
 \frac{\,1\,}{t^{\nu/\beta}}
 \Phi_2\left(\frac{\,d(x,y)\,}{t^{1/\beta}} \right).
\end{align*}
Then we have the following:
\begin{enumerate}
\item Suppose $\mu(E)<+\infty$. Then
$\mu\in S_K^{\,p,\delta}$
for any $\delta\in]0,(\eta-p(\nu-\beta))/p\beta[$.
\item Suppose $p>1$ and $\mu\in S_D$.  
Then $\mu\in S_K^{\,p,\delta}$ holds for 
any $\delta\in]0,(\eta-p(\nu-\beta))/p\beta[$ 
{\rm(}resp.~$\delta\in]0,(p-1)(\eta-p(\nu-\beta))/p^2\beta[${\rm)}
under $p\in [(\eta+\beta)/\nu,\eta/(\nu-\beta)_+[$ {\rm(}resp.~$p\in ]1,(\eta+\beta)/\nu[${\rm)}. Here 
$\eta/(\nu-\beta)_+:=\eta/(\nu-\beta)$ for $\nu>\beta$, and 
$\eta/(\nu-\beta)_+:=+\infty$ for $\nu\leq\beta$. 
%If $\nu\geq\eta+\beta$, then $\mu\in S_K^{\,p,\delta}$
%for any $\delta\in]0,(\eta-p(\nu-\beta))/p\beta[$ under $p\in]1,\eta/(\nu-\beta)[$.
\item 
Suppose $p>1$ and $\mu\in S_D$. 
Assume further that
there exists $C>0$ such that
$\sup_{x\in E}\mu(B_r(x))\leq Cr^{\eta}$ holds
for any $r\in]0,+\infty[$, 
or
$
 \int_1^\infty
  u^{\nu+\gamma-1}\Phi_2(u)
 \d u
<
 +\infty
$
holds for any $\gamma>0$. 
Then $\mu\in S_K^{\,p,\delta}$
for any $\delta\in]0,(\eta-p(\nu-\beta))/p\beta[$.
\end{enumerate}
\end{thm}

\section{Proofs of theorems and corollaries}\label{sec:proof} 

\subsection{Proof of Theorem \ref{Th1}}
To prove Theorem \ref{Th1}, we begin with auxiliary lemmas.
\begin{lem}[{\cite[Lemma~4.1]{KwT:Katounderheat}}]\label{lem:lowerestimate}
Under {\bf (A)}\ref{asmp:Phi}, there exists $C_{\nu,\beta,t_0}^{'}>0$ such that for any $t\in]0,t_0]$ 
{\rm(}$t\in]0,+\infty[$ if $t_0=+\infty${\rm)}  
\begin{enumerate}
\item\label{item:lower<} for $\nu < \beta$ and $x,y\in E$ with $d(x,y)^{\beta}< t$, we have   
\begin{equation*}
 \int^{t}_{0} p_s(x,y) \d s
\geq
 C_{\nu,\beta,t_0}^{'}
 t^{1-\frac{\,\nu\,}{\beta}},
\end{equation*}
\item\label{item:lower=} for $\nu=\beta$ and $x,y\in E$ with $d(x,y)^{\beta/2}< t < 1/2$, 
 we have 
\begin{equation*}
 \int^{t}_{0} p_s(x,y) \d s
\geq
 C_{\nu,\beta,t_0}^{'}
 \log(d(x,y)^{-1}),
\end{equation*}
\item\label{item:lower>} for $\nu > \beta$ and $x,y\in E$ with $d(x,y)^{\beta}< t$,  
\begin{equation*}
 \int^{t}_{0} p_s(x,y) \d s
\geq
 C_{\nu,\beta,t_0}^{'}d(x,y)^{\beta-\nu}.
\end{equation*}
\end{enumerate}
\end{lem}

\begin{lem}[{\cite[Lemma~4.2]{KwT:Katounderheat}}]\label{lem:lowerestimatelargetime}
Under {\bf (A)}\ref{asmp:Phi}, for any $t\in]0,+\infty[$, 
there exists $C_{\nu,\beta,t_0,t}^{'}>0$ such that 
for $x,y\in E$ with $d(x,y)^{\beta}< t$, we have 
\begin{equation*}
 \int^{t}_{0} p_s(x,y) \d s \geq  C_{\nu,\beta,t_0,t}^{'}.
\end{equation*}
\end{lem}

\begin{lem}[{\cite[Lemma~4.3]{KwT:Katounderheat}}]\label{lem:upperestimate}
Under {\bf (A)}\ref{asmp:Phi}, there exists $C_{\nu,\beta,t_0}>0$ such that for any $t\in]0,t_0]$ 
{\rm(}$t\in]0,+\infty[$ if $t_0=+\infty${\rm)}  
\begin{enumerate}
\item\label{item:upper<} 
for $\nu<\beta$ and $x,y\in E$, we have  
\begin{equation*}
 \int^{t}_{0} p_s(x,y) \d s
\leq
 C_{\nu,\beta,t_0}
 t^{1-\nu/\beta},
\end{equation*}
\item\label{item:upper=} for $\nu=\beta$ and $x,y\in E$ with $d(x,y)^{\beta}\lor t < 1/2$, we have  
\begin{equation*}
 \int^{t}_{0} p_s(x,y) \d s
\leq
 C_{\nu,\beta,t_0}
 \log(d(x,y)^{-1}),
\end{equation*}
\item\label{item:upper>} for $\nu>\beta$ and $x,y\in E$, we have 
\begin{equation*}
 \int^{t}_{0} p_s(x,y) \d s
\leq
 C_{\nu,\beta,t_0}
 d(x,y)^{\beta-\nu}.
\end{equation*}
\end{enumerate}
\end{lem}

\begin{apf}{Theorem~\ref{Th1}} 
The implications
\ref{item:resolventany}$\Longrightarrow$%
\ref{item:resolventsome}$\Longrightarrow$%
\ref{item:heatany}$\Longrightarrow$%
\ref{item:heatsome}
and 
$(2')$$\Longrightarrow$%
$(3')$$\Longrightarrow$%
$(4')$$\Longrightarrow$%
$(5')$
are trivial in view of the estimate
\begin{align*}
\int_0^tp_s(x,y)\d s\leq e^{\alpha t}r_{\alpha}(x,y).
\end{align*} 
First
we show \ref{item:heatsome}$\Longrightarrow$\ref{item:Green}. 
Suppose \ref{item:heatsome}. Then 
\begin{align*}
\lim_{r\to0}\sup_{x\in E}\int_{B_r(x)}
\left(\int_0^{t_1}
p_s(x,y)\d s\right)^p\mu(\d y)=0
\end{align*}
holds for some $t_1>0$. 
We may assume $t_0<1/2\land t_1$. By Lemma~\ref{lem:lowerestimate}, 
we see that for $d(x,y)<r<t_0^{2/\beta}/2$
\begin{align*}
 C_{\nu,\beta,t_0}'G(x,y)
\leq
 \int_0^{t_0}p_s(x,y)\d s
\leq
 \int_0^{t_1}p_s(x,y)\d s
.
\end{align*}
Then 
\begin{align*}
\left(C_{\nu,\beta,t_0}'\right)^p\lim_{r\to0}\sup_{x\in E}\int_{B_r(x)}G(x,y)^p\mu(\d y)\leq
\lim_{r\to0}\sup_{x\in E}\int_{B_r(x)}
\left(\int_0^{t_1}
p_s(x,y)\d s\right)^p\mu(\d y)=0.
\end{align*}
Thus we have \ref{item:Green}.
The proof of 
$(5')$$\Longrightarrow$$(1')$ is similar.  
Next we 
show \ref{item:Green}$\Longrightarrow$\ref{item:resolventany}. 
Owing to the estimate in the proof of  Lemma~\ref{lem:upperestimate}\ref{item:upper>}, 
for $\nu>\beta$ 
we have  
\begin{align*}
 r_{\alpha}(x,y)
&=
 \sum_{k=0}^{\infty}
 \int_{kt_0}^{(k+1)t_0}e^{-\alpha s}p_s(x,y)\d s
\\
&\leq
 \sum_{k=0}^{\infty}
 e^{-\alpha kt_0}
 \int_0^{t_0}p_{s+kt_0}(x,y)\d s
\\
&=
 \sum_{k=0}^{\infty}
 e^{-\alpha kt_0}
 \int_0^{t_0}\int_E p_{kt_0}(x,z)p_s(z,y)\m(\d z)\d s
\\
&\leq
 \sum_{k=0}^{\infty}
 e^{-\alpha kt_0}
 \int_E p_{kt_0}(x,z)
  \frac{\,M_{\nu,\beta}\,}{d(z,y)^{\nu-\beta}}
 \m(\d z),
\end{align*}
where $M_{\nu,\beta}:=\beta\int_0^{\infty}u^{\nu-\beta-1}\Phi_2(u)\d u<+\infty$. 
Hence 
\begin{align*}
\left(\int_{B_r(x)}r_{\alpha}(x,y)^p\mu(\d y) \right)^{\frac{\,1\,}{p}}
&\leq \sum_{k=0}^{\infty}e^{-\alpha kt_0}M_{\nu,\beta}
\left(\int_{B_r(x)}\left(\int_Ep_{kt_0}(x,z)\cdot\frac{\,\m(\d z)\,}{d(z,y)^{\nu-\beta}} \right)^p\mu(\d y) \right)^{\frac{\,1\,}{p}}\\
&
\leq \sum_{k=0}^{\infty}e^{-\alpha kt_0}M_{\nu,\beta}
\left\{
\left(\int_{B_r(x)}\left(\int_{B_{2r}(x)}p_{kt_0}(x,z)\cdot\frac{\,\m(\d z)\,}{d(z,y)^{\nu-\beta}} \right)^p\mu(\d y) \right)^{\frac{\,1\,}{p}}\right. \\
&\hspace{1cm}\left.
+
\left(\int_{B_r(x)}\left(\int_{B_{2r}(x)^c}p_{kt_0}(x,z)\cdot\frac{\,\m(\d z)\,}{d(z,y)^{\nu-\beta}} \right)^p\mu(\d y) \right)^{\frac{\,1\,}{p}}
 \right\}.
\end{align*} 
The first term in the curly brackets of the right-hand side
is bounded from above by
\begin{align*}
 \left(
  \int_{B_r(x)}
  \int_{B_{2r}(x)}
   p_{kt_0}(x,z)\frac{\,\m(\d z)\,}{d(z,y)^{p(\nu-\beta)}}
  \mu(\d y)
 \right)^{\frac{\,1\,}{p}} 
\leq&
 \left(
  \int_{B_{2r}(x)}
  \int_{B_{3r}(z)}
   p_{kt_0}(x,z)
   \frac{\,\mu(\d y)\,}{d(z,y)^{p(\nu-\beta)}}
  \m(\d z)
 \right)^{\frac{\,1\,}{p}} 
\\
\leq&
 \left(
  \sup_{z\in E} 
  \int_{B_{3r}(z)}\frac{\,\mu(\d y)\,}{d(z,y)^{p(\nu-\beta)}}
 \right)^{\frac{\,1\,}{p}} 
\end{align*}
and the second term is bounded from above by
\begin{align*}
 \left(
  \int_{B_r(x)}
   \frac{\,1\,}{r^{p(\nu-\beta)}}
   \left(
    \int_{B_{r}(y)^c}
     p_{kt_0}(x,z)
    \m(\d z)
   \right)^p
  \mu(\d y)
 \right)^{\frac{\,1\,}{p}}
&\leq
 \left(
  \int_{B_r(x)}
  \frac{\,\mu(\d y)\,}{d(x,y)^{p(\nu-\beta)}}
 \right)^{\frac{\,1\,}{p}}.
\end{align*}
Thus, for $\nu>\beta$ 
\begin{align}
\label{eq:1->2_bound_nu>beta}
 \int_{B_r(x)}r_{\alpha}(x,y)^p\mu(\d y)
\leq
 \left(
  \frac{\,2M_{\nu,\beta}\,}{1-e^{-\alpha t_0}}
 \right)^p
 \sup_{x\in E}
 \int_{B_{3r}(x)}G(x,y)^p\mu(\d y). 
\end{align}
Suppose $\nu=\beta$.
We can see that for $d(y,z)\geq t_0^{1/\beta}$,
\begin{align*}
 \beta\int_{d(z,y)/t_0^{1/\beta}}^{\infty}u^{-1}\Phi_2(u)\d u
\leq
 \beta\int_{1}^{\infty}u^{-1}\Phi_2(u)\d u
\leq
 \beta\int_{1}^{\infty}u^{\nu-1}\Phi_2(u)\d u
\end{align*}
and for $d(y,z)< t_0^{1/\beta} (<1)$,
\begin{align*}
 \beta\int_{d(z,y)/t_0^{1/\beta}}^{\infty}u^{-1}\Phi_2(u)\d u
=&
 \beta\int_{d(z,y)/t_0^{1/\beta}}^{1}u^{-1}\Phi_2(u)\d u
+
 \beta\int_{1}^{\infty}u^{-1}\Phi_2(u)\d u
\\
\leq&
 \beta\Phi_2(0)
 \log(d(y,z)^{-1})
+
 \beta\int_{1}^{\infty}u^{\nu-1}\Phi_2(u)\d u.
\end{align*}
Then we have that
\begin{align*}
 r_{\alpha}(x,y)
&=
 \sum_{k=0}^{\infty}
 e^{-\alpha k t_0}
 \int_E
  p_{kt_0}(x,z)
  \left(\int_0^{t_0}p_s(z,y)\d s\right)
 \m(\d z)
\\
&\leq
 \sum_{k=0}^{\infty}
 e^{-\alpha k t_0}
 \int_E
  p_{kt_0}(x,z)
  \left(\beta\int_{d(z,y)/t_0^{1/\beta}}^{\infty}u^{-1}\Phi_2(u)\d u\right)
 \m(\d z)
\\ 
%&=
% \sum_{k=0}^{\infty}
% e^{-\alpha k t_0}
% \int_{d(y,z)\geq t_0^{1/\beta}}
%  p_{kt_0}(x,z)
%  \left(\beta\int_{d(z,y)/t_0^{1/\beta}}^{\infty}u^{-1}\Phi_2(u)\d u\right)
% \m(\d z)
%\\ 
%&+
% \sum_{k=0}^{\infty}
% e^{-\alpha k t_0}
% \int_{d(y,z)< t_0^{1/\beta}}
%  p_{kt_0}(x,z)
%  \left(\beta\int_{d(z,y)/t_0^{1/\beta}}^{\infty}u^{-1}\Phi_2(u)\d u\right)
% \m(\d z)
%\\ 
&\leq
 \sum_{k=0}^{\infty}
 e^{-\alpha k t_0}
 \left(
  2
  \beta\int_1^{\infty}u^{\nu-1}\Phi_2(u)\d u
 \right.
\\ 
&\hspace{1cm}
 \left.
 +
  \beta
  \Phi_2(0)
  \int_{d(y,z)<t_0^{1/\beta}}p_{kt_0}(x,z)
  \log(d(y,z)^{-1})\m(\d z)
 \right).
\end{align*}
In the same way
to obtain \eqref{eq:1->2_bound_nu>beta},
we have that 
\begin{align*}
&
 \int_{B_r(x)}
 \left(
   \beta
   \Phi_2(0)
   \int_{d(y,z)<t_0^{1/\beta}}p_{kt_0}(x,z)
   \log(d(y,z)^{-1})\m(\d z)
 \right)^p
 \mu(\d y)
\\
\leq&
 \left(2 \beta \Phi_2(0)\right)^p
 \sup_{x\in E}\int_{B_{3r}(x)}G(x,y)^p\mu(\d y). 
\end{align*}
Hence,
for $\nu=\beta$ with $r<e^{-1}$
\begin{align*}
 \int_{B_r(x)}r_{\alpha}(x,y)^p\mu(\d y)
\leq
 \left(
  \frac{\, 2C_1+2C_2\,}{1-e^{-\alpha t_0}}
 \right)^p
 \sup_{x\in E}
 \int_{B_{3r}(x)}G(x,y)^p\mu(\d y),  
\end{align*}
where
$C_1:=\beta\int_1^{\infty}u^{\nu-1}\Phi_2(u)\d u$
and $C_2:=\beta\Phi_2(0)$. 
Here we use that for $e<r^{-1}$
\begin{align*}
\mu(B_r(x))\leq (\log r^{-1})^p\mu(B_r(x))\leq\int_{B_r(x)}(\log d(x,y)^{-1})^p\mu(\d y).
\end{align*} 
Therefore we obtain that for a constant $D_{\nu,\beta,\alpha,t_0}>0$ 
\begin{align*}
\lim_{r\to0}\sup_{x\in E}\int_{B_r(x)}r_{\alpha}(x,y)^p\mu(\d y)\leq D_{\nu,\beta,\alpha,t_0}\lim_{r\to0}
\sup_{x\in E}\int_{B_{r}(x)}G(x,y)^p\mu(\d y)=0,  
\end{align*}
which implies the desired assertion. 
The proof of $(1')$$\Longrightarrow$$(2')$ is also similar.  
\end{apf}

\subsection{Proof of Theorem \ref{Th2}}

To prove Theorem \ref{Th2}, we begin with auxiliary lemmas.
\begin{lem}\label{lem:lowerPhi}
Under {\bf (A)}\ref{asmp:Phi}, we have 
$S_D^{\,p}\subset D_{\nu,\beta}^{\,p}$ and 
$S_K^{\,p}\subset K_{\nu,\beta}^{\,p}$.   
Moreover, $\mu\in S_D^{\,p}$ implies $\sup\limits_{x\in E}\mu(B_R(x))<+\infty$ for all $R>0$. 
In particular, $S_D^{\,p}\subset K_{\nu,\beta}^{\,p}$ if $\nu<\beta$. 
\end{lem}

\begin{pf} 
Take $\mu\in S_D^{\,p}$. 
By Lemma~\ref{lem:lowerestimate}, for $\nu > \beta$ (resp.~$\nu=\beta$) with $r:=t^{1/\beta}$ 
(resp.~$r:=t^{2/\beta}$), we have  
\begin{equation*}
 \int_E
  \left(
   \int_0^t p_s(x,y)\d s
  \right)^p
 \mu(\d y)
\geq
 \left(C_{\nu,\beta,t_0}^{'}\right)^p
 \int_{d(x,y)< r} G(x,y)^p\mu(\d y).
\end{equation*}
Then we see $\mu\in D_{\nu,\beta}^{\,p}$. 
Hence, for $\mu\in S_K^{\,p}$, we have 
\begin{align*}
\lim_{r\downarrow 0} \sup_{x\in E} \int_{d(x,y)< r} G(x,y)^p\mu(\d y) 
&\leq \frac{\,1\,}{\left(C_{\nu,\beta,t_0}^{'}\right)^p} \lim_{t\downarrow 0} \sup_{x\in E} 
\int_E\left(\int_0^t p_s(x,y)\d s\right)^p\mu(\d y)=0. 
\end{align*} 
Take $\mu\in S_D^{\,p}$. 
By using Lemma~\ref{lem:lowerestimatelargetime}, 
\begin{align*}
\int_E \left(\int_0^t p_s(x,y)\d s\right)^p\mu(\d y) &\geq 
\int_{d(x,y)< t^{1/\beta}}\left(\int_0^t p_s(x,y)\d s\right)^p\mu(\d y) \geq 
\left(C_{\nu,\beta,t_0,t}^{'}\right)^p\int_{d(x,y)< t^{1/\beta}}\mu(\d y).
\end{align*}
So it suffices to apply Lemma~\ref{lem:Dynkinequi} with $t=r^{\beta}$. 
\end{pf}

\begin{lem}\label{lem:upperPhi}
Under {\bf (A)}\ref{asmp:Phi}, we have $D_{\nu,\beta}^{\,p}\subset S_1$. 
\end{lem}

\begin{pf}
Recall that every $\mu\in D_{\nu,\beta}^{\,p}$ is a positive Radon measure on $E$, that is, $\mu(K)<+\infty$ 
for each compact set $K$. It suffices to show that for each compact set $K$, 
$\mu\in D_{\nu,\beta}^{\,p}$ implies $\1_K\mu\in S_D^{\,p}$. 
In fact, $\1_K\mu\in S_D^{\,p}$ implies $\1_K\mu\in S_{00}$, hence $\mu\in S_1$ 
by \cite[Theorem~5.1.7(iii)]{FOT}.

By Lemma~\ref{lem:upperestimate}\ref{item:upper<}, for $\nu < \beta$, we have 
\begin{equation*}
 \int_K
  \left(
   \int_0^t p_s(x,y)\d s
  \right)^p
 \mu(\d y)
\leq
 \left(C_{\nu,\beta,t_0}\right)^p
 \mu(K)
 t^{p\left(1-\nu/\beta\right)}.
\end{equation*}
Then we obtain $\1_K\mu\in S_D^{\,p}$ in this case. 
For $\nu \geq  \beta$ with $t=r^{\beta}<1/2$ and $r<1/e$, 
we have from Lemma~\ref{lem:upperestimate}\ref{item:upper=},\ref{item:upper>},                                   
\begin{equation*}
 \int_{K\cap B_r(x)}
  \left(
   \int_0^t p_s(x,y)\d s
  \right)^p
 \mu(\d y)
\leq
 \left(C_{\nu,\beta,t_0}\right)^p
 \int_{d(x,y)< r}G(x,y)^p\mu(\d y)
\end{equation*}
and  
\begin{align*} 
\int_{K\cap B_r(x)^c}\left(\int_0^t p_s(x,y)\d s\right)^p \mu(\d y)
&\leq 
\left(\beta r^{\beta-\nu} \int_{r/t^{1/\beta}}^{\infty}u^{\nu-1}\Phi_2(u)\d u \right)^p \mu(K) \\
&=
\left(\beta r^{\beta-\nu} \int_{1}^{\infty}u^{\nu-1}\Phi_2(u)\d u \right)^p \mu(K).
\end{align*}
Thus there exists  $t\in]0,t_0\land\frac12[$ such that 
\begin{equation*}
 \sup_{x\in E}
 \int_E
  \left(\int_0^t p_s(x,y)\d s\right)^p \1_K(y)
 \mu(\d y)
<+\infty,
\end{equation*}
which implies $\1_K\mu\in S_D^{\,p}$. 
\end{pf}

\begin{apf}{Theorem~\ref{Th2}}
The assertion $K_{\nu,\beta}^{\,p}=S_K^{\,p}$ 
of Theorem~\ref{Th2} is proved under $p=1$ 
by \cite[Theorem~3.2]{KwT:Katounderheat}. That is, 
we already know $K_{\nu,\beta}=S_K\subset S_D$ under 
{\bf (A)}\ref{asmp:SC}, {\bf (A)}\ref{asmp:Bishop} and {\bf (A)}\ref{asmp:Phi}. 
We do not prove the coincidence $D_{\nu,\beta}=S_D$ under 
{\bf (A)}\ref{asmp:SC}, {\bf (A)}\ref{asmp:Bishop} and {\bf (A)}\ref{asmp:Phi} in  \cite[Theorem~3.2]{KwT:Katounderheat}. But the method of the proof of 
\cite[Theorem~3.2]{KwT:Katounderheat} still works to 
prove $D_{\nu,\beta}=S_D$. We omit the details but note that its proof can be achieved by use of the estimates  
\begin{align*}
\sup_{x\in E}{\EE}_x[A_{\tau_{B_r(x)}\land s}^{\mu}]\leq 
\left\{\begin{array}{cc} \displaystyle{ C_{\nu,\beta,t_0}\sup_{x\in E}\int_{B_r(x)}G(x,y)\mu(\d y)}, & \nu\geq\beta, \\ \displaystyle{
C_{\nu,\beta,t_0}s^{1-\nu/\beta}\sup_{x\in E}\mu(B_1(x))}, & \nu<\beta\end{array}\right.
\end{align*}
for $s\in]0,t_0\land 1/2[$
and 
\begin{align*}
\sup_{x\in E}{\PP}_x(T_k<s_0)<(\gamma+\eps)^k.
\end{align*}
Here $A^{\mu}$ is the positive 
continuous additive functional associated to $\mu$, 
$\gamma$ is the constant appeared in {\bf (A)}\ref{asmp:SC} 
satisfying $\gamma+\eps<1$ for some $\eps>0$, and 
$(T_k)_{k\geq0}$ is the sequence of stopping times
defined by $T_0=0$, 
$T_{k+1}=T_k+(\tau_{B_{r_0}(X_0)}\land s)\circ \theta_{T_k}$
(see \cite[pp.~102]{KwT:Katounderheat}).

By Lemma~\ref{lem:monotonicityKato}, 
we have 
$D_{\nu,\beta}^{\,p}\subset D_{\nu,\beta}=S_D$ 
and 
$K_{\nu,\beta}^{\,p}\subset K_{\nu,\beta}=S_K$. 
From Lemma~\ref{lem:lowerPhi}, 
we also have $S_D^{\,p}\subset D_{\nu,\beta}^{\,p}$ 
and 
$S_K^{\,p}\subset K_{\nu,\beta}^{\,p}$.  
So it suffices to show
$D_{\nu,\beta}^{\,p}\subset S_D^{\,p}$ and 
$K_{\nu,\beta}^{\,p}\subset S_K^{\,p}$.  
Take $\mu\in K_{\nu,\beta}^{\,p}$ and fix $r>0$.
Note that for $d(x,y)\geq r^{1/\beta}$, 
\begin{align*}
 \int_0 ^t
  p_s(x, y)
 \d s
\leq
 \int_{(r/t)^{1/\beta}}^{\infty}
  u^{-1}\Phi_2(u)
 \d u
\longrightarrow 0 \quad \text{ as }\quad t\to0,
\end{align*}
because $u^{-1}\Phi_2(u)$ is integrable on $[1,+\infty[$
under the condition {\bf (A)}\ref{asmp:Bishop}.
By combining this with the fact
$\mu\in K_{\nu,\beta}^{\,p} \subset S_D$ as noted above,
we have 
\begin{align*}
 \sup_{x\in E}&
 \int_{d(x, y)\geq r^{1/\beta}}
  \left(
   \int_0 ^t
    p_s(x, y)
   \d s
  \right)^p
 \mu(\d y)
\\
&\leq
 \left(
  \int_{(r/t)^{1/\beta}}^{\infty}
   u^{-1}
   \Phi_2(u)
  \d u
 \right)^{p-1}
 \sup_{x\in E}
 \int_E
  \left(
   \int_0 ^t
    p_s(x, y)
   \d s
  \right)
 \mu(\d y)
\longrightarrow0\quad\text{ as }\quad t\to0. 
\end{align*}
Hence,
\begin{align*}
 \varlimsup_{t \to 0}&
 \sup_{x\in E}
 \int_{E}
  \left(
   \int_0 ^t
    p_s(x, y)
   \d s
  \right)^p
 \mu(\d y)
\\
&\leq
 \varlimsup_{t \to 0}
 \sup_{x\in E}
 \int_{d(x, y)< r^{1/\beta}}
  \left(
   \int_0 ^t
    p_s(x, y)
   \d s
  \right)^p
 \mu(\d y)
+
 \varlimsup_{t \to 0}
 \sup_{x\in E}
 \int_{d(x, y)\geq r^{1/\beta}}
  \left(
   \int_0 ^t
    p_s(x, y)
   \d s
  \right)^p
 \mu(\d y)
\\
&\leq
 \sup_{x\in E}
 \int_{d(x, y)< r^{1/\beta}}
  \left(
   \int_0 ^1
    p_s(x, y)
   \d s
  \right)^p
 \mu(\d y),
\end{align*}
and the right-hand side of the above inequality goes to zero
as $r\rightarrow 0$ by Theorem~\ref{Th1},
which concludes $\mu\in S_K^{\,p}$. 
The proof of
$\mu\in  D_{\nu,\beta}\subset S_D$ $\Longrightarrow$
$\mu\in S_D^{\,p}$ is similar.

Therefore we obtain 
$S_D^{\,p}=D_{\nu,\beta}^{\,p}$ and 
$S_K^{\,p}=K_{\nu,\beta}^{\,p}$ under {\bf (A)}\ref{asmp:SC}, {\bf (A)}\ref{asmp:Bishop} and {\bf (A)}\ref{asmp:Phi}. 
Finally, we prove the rest assertions. 
As proved in 
Lemma~\ref{lem:lowerPhi} 
(see \cite[Lemma~4.4]{KwT:Katounderheat}), 
every $\mu\in S_D^{\,p}$ satisfies 
\eqref{eq:ballunifinite} 
under 
 {\bf (A)}\ref{asmp:Phi}. The same holds for $\mu\in K_{\nu,\beta}^{\,p}$ under {\bf (A)}\ref{asmp:SC}, {\bf (A)}\ref{asmp:Bishop} and {\bf (A)}\ref{asmp:Phi}. 
 When $\nu<\beta$, we have $K_{\nu,\beta}^{\,p}=K_{\nu,\beta}=S_D=S_K=S_K^{\,p}=S_D^{\,p}$ in view of 
Lemma~\ref{lem:lowerPhi} under the above assumptions. In this case, 
any Radon measure $\mu$ satisfying \eqref{eq:ballunifinite} always belongs to $K_{\nu,\beta}^{\,p}$. 
\end{apf}

\subsection{Proofs of Theorem \ref{Th3}, Corollary \ref{cor:Th3} and Theorem \ref{Th4}}

\begin{apf}{Theorem~\ref{Th4}} 
Suppose $\nu-p(\nu-\beta)>0$. 
When $\nu\geq\beta$,
take $\gamma>0$ such that
$\alpha/q<\gamma<\nu-p(\nu-\beta)$.
Then the conclusion immediately follows from
Lemmas \ref{lem:Maq_sub_Mg1} and \ref{lem:SchechtersubKato}
stated below.

It remains to show the assertion in the case $\nu<\beta$.
First
we assume $\nu\geq\alpha$.
When
$q\in ]1, +\infty[$, H\"older's inequality gives that
\begin{align*}
 \sup_{x\in E}&
 \int_{d(x,y)\leq1}|f(y)|\m(\d y)
\\
\leq&
 \left(
  \sup_{x\in E}
  \int_{d(x,y)\leq1}
   d(x,y)^{\frac{\,\nu-\alpha\,}{q-1}}
  \m(\d y)
 \right)^{\frac{\,q-1\,}{q}}
 \left(
  \sup_{x\in E}
  \int_{d(x,y)\leq1}
   \frac{\,|f(y)|^q\,}{d(x,y)^{\nu-\alpha}}
  \m(\d y)
 \right)^{\frac{\,1\,}{q}}
\\
\leq&
 V(2)^{\frac{\,q-1\,}{q}}
 \left(
  \sup_{x\in E}
  \int_{d(x,y)\leq1}
   \frac{\,|f(y)|^q\,}{d(x,y)^{\nu-\alpha}}
  \m(\d y)
 \right)^{\frac{\,1\,}{q}}
<
 +\infty
. 
\end{align*}
When $q=1$, we have
\begin{align*}
 \sup_{x\in E}
 \int_{d(x,y)\leq1}|f(y)|\m(\d y)
=&
 \sup_{x\in E}
 \int_{d(x,y)\leq1}
  \frac{\,|f(y)|\,}{d(x,y)^{\nu-\alpha}}d(x,y)^{\nu-\alpha}
 \m(\d y)
\\
\leq&
 \sup_{x\in E}
 \int_{d(x,y)\leq1}
  \frac{\,|f(y)|\,}{d(x,y)^{\nu-\alpha}}
 \m(\d y)
<
 +\infty
.
\end{align*}
Next we assume $\nu<\alpha$. 
For $q\in]\alpha/\nu,+\infty[$,
H\"older's inequality gives that
%\begin{align*}
% \sup_{x\in E}&
% \int_{d(x,y)\leq1}
%  |f(y)|
% \m(\d y)
%\\
%\leq&
% \left(
%  \sup_{x\in E}
%  \int_{d(x,y)\leq1}d(x,y)^{\frac{\,\nu-\alpha\,}{q-1}} \m(\d y)
% \right)^{\frac{\,q-1\,}{q}}
% \left(
%  \sup_{x\in E}
%  \int_{d(x,y)\leq1}
%   \frac{\,|f(y)|^q\,}{d(x,y)^{\nu-\alpha}}
%  \m(\d y)
% \right)^{\frac{\,1\,}{q}}
%\\
%\leq&
% \left(
%  \sup_{x\in E}
%  \sum_{k=0}^{\infty}
%  \int_{\frac{\,1\,}{2^{k+1}}\leq d(x,y)<\frac{\,1\,}{2^k}}
%  d(x,y)^{\frac{\,\nu-\alpha\,}{q-1}}\m(\d y)
% \right)^{\frac{\,q-1\,}{q}}
% \left(
%  \sup_{x\in E}
%  \int_{d(x,y)\leq1}
%   \frac{\,|f(y)|^q\,}{d(x,y)^{\nu-\alpha}}
%  \m(\d y)
% \right)^{\frac{\,1\,}{q}}
%\\
%\leq&
% \left(
%  C
%  \sum_{k=0}^{\infty}
%  (2^{k+1})^{\frac{\,\alpha-\nu\,}{q-1}}2^{-k\nu}
% \right)^{\frac{\,q-1\,}{q}}
% \left(
%  \sup_{x\in E}
%  \int_{d(x,y)\leq1}
%   \frac{\,|f(y)|^q\,}{d(x,y)^{\nu-\alpha}}
%  \m(\d y)
% \right)^{\frac{\,1\,}{q}}
%\\
%=&
% \left(
%  C\cdot 
%  2^{\frac{\,\alpha-\nu\,}{q-1}}
%  \sum_{k=0}^{\infty}
%  2^{-k\frac{\,\nu q-\alpha\,}{q-1}}
% \right)^{\frac{\,q-1\,}{q}}
% \left(
%  \sup_{x\in E}
%  \int_{d(x,y)\leq1}
%   \frac{\,|f(y)|^q\,}{d(x,y)^{\nu-\alpha}}
%  \m(\d y)
% \right)^{\frac{\,1\,}{q}}
%<
% +\infty
%.
%\end{align*}
\begin{align*}
 \sup_{x\in E}&
 \int_{d(x,y)\leq1}
  |f(y)|
 \m(\d y)
\\
\leq&
 \left(
  \sup_{x\in E}
  \int_{d(x,y)\leq1}
   d(x,y)^{\frac{\,\nu-\alpha\,}{q-1}}
  \m(\d y)
 \right)^{\frac{\,q-1\,}{q}}
 \left(
  \sup_{x\in E}
  \int_{d(x,y)\leq1}
   \frac{\,|f(y)|^q\,}{d(x,y)^{\nu-\alpha}}
  \m(\d y)
 \right)^{\frac{\,1\,}{q}}
.
\end{align*}
The quantity in the first parentheses 
of the right-hand side
is bounded from above by
\begin{align*}
 \sup_{x\in E}&
 \sum_{k=0}^{\infty}
 \int_{\frac{\,1\,}{2^{k+1}}\leq d(x,y)<\frac{\,1\,}{2^k}}
  d(x,y)^{\frac{\,\nu-\alpha\,}{q-1}}
 \m(\d y)
\\
&\leq
 C
 \sum_{k=0}^{\infty}
 (2^{k+1})^{\frac{\,\alpha-\nu\,}{q-1}}2^{-k\nu}
=
 C\cdot 
 2^{\frac{\,\alpha-\nu\,}{q-1}}
 \sum_{k=0}^{\infty}
 2^{-k\frac{\,\nu q-\alpha\,}{q-1}}
<
 +\infty.
\end{align*}
Hence $|f|\d\m\in K_{\nu,\beta}^{\,p}$.
\end{apf}

The following lemmas are extensions of 
\cite[Proposition 4.1]{AizSim:Brown}
and the inclusion below 
\cite[Proposition 4.1]{AizSim:Brown},
in which the results are obtained for the case of $p=1$ 
and for a Brownian motion on $\mathbb{R}^d$ with $d\geq 3$:
\begin{lem}
\label{lem:Maq_sub_Mg1}
Assume {\bf (A)}\ref{asmp:Bishop} and let
$\alpha, \gamma >0$ and $q\geq 1$ with $\alpha < \gamma q$.
Then it holds that
$M_{\alpha, q}\subset M_{\gamma, 1}$.
\end{lem}

\begin{pf}
The conclusion is trivial when $q=1$.
Suppose $q>1$ and set
$
 a
:=
 \nu
+
 \frac{\,q\,}{q-1}
 \bigl(
  \gamma - \nu + \frac{\,\nu-\alpha\,}{q}
 \bigr)
=
 \frac{\,1\,}{q-1} (q\gamma -\alpha)
>
 0
$.
Consider the integral
$
 \int_{d(x, y)\leq 1}
  d(x, y)^{a-\nu}
 \m(\d y)
$.
Note that,
by {\bf (A)}\ref{asmp:Bishop} there exists $C>0$ such that
$\sup_{x\in E}\m(B_r(x))\leq V(r)\leq Cr^{\nu}$
for any $r\in]0,2]$.
When $a\geq \nu$, the integral is bounded from above by
$\m(B_1(x)) \leq V(2) < +\infty$.
When $a<\nu$, the integral is bounded from above by
\begin{align*}
 \sum_{k=0}^\infty
 \int_{2^{-(k+1)}<d(x, y)\leq 2^{-k}}
  d(x, y)^{a-\nu}
 \m(\d y)
&\leq
 \sum_{k=0}^\infty
 2^{-(k+1)(a-\nu)} \m(B_{2^{-k}}(x))
\\
&\leq
 C
 2^{-(a-\nu)}
 \sum_{k=0}^\infty
 2^{-ak}
<
 +\infty.
\end{align*}
Now, let $f\in M_{\alpha, q}$.
By H\"older's inequality we have
\begin{align*}
 \int_{d(x, y)\leq 1}&
  \frac{\,|f(y)|\,}{d(x, y)^{\nu-\gamma}}
 \m(\d y)
\\
&=
 \int_{d(x, y)\leq 1}
  \frac{\,|f(y)|\,}{d(x, y)^{\frac{\,\nu-\alpha\,}{q}} }
  d(x, y)^{\gamma-\nu+\frac{\,\nu-\alpha\,}{q}}
 \m(\d y)
\\
&\leq
 \left(
  \int_{d(x, y)\leq 1}
   \frac{\,|f(y)|^q\,}{d(x, y)^{\nu-\alpha}}
  \m(\d y)
 \right)^{\frac{\,1\,}{q}}
 \left(
  \int_{d(x, y)\leq 1}
   d(x, y)^{a-\nu}
  \m(\d y)
 \right)^{\frac{\,q-1\,}{q}},
\end{align*}
which concludes $f\in M_{\gamma, 1}$.
\end{pf}

\begin{lem}
\label{lem:SchechtersubKato}
Let $p\in[1,+\infty[$ and $\nu\geq\beta$ with $\nu-p(\nu-\beta)>0$.
For any $f \in  M_{\gamma,1}$ $(\gamma>0)$, 
we have $|f|\d\m\in K_{\nu,\beta}^{\,p}$ if $\gamma < \nu-p(\nu-\beta)$.
\end{lem}

\begin{pf}
Take $f\in M_{\gamma, 1}$.
When $\nu>\beta$, we have for $r<1$
\begin{align*}
 \int_{B_r(x)}
  \frac{\,|f(y)|\,}{d(x, y)^{p(\nu-\beta)}}
 \m(\d y)
&=
 \int_{B_r(x)}
  \frac{\,|f(y)|\,}{d(x, y)^{\nu-\gamma}}
  d(x, y)^{\nu-\gamma-p(\nu-\beta)}
 \m(\d y)
\\
&\leq
 r^{\nu-\gamma-p(\nu-\beta)}
 \int_{B_r(x)}
  \frac{\,|f(y)|\,}{d(x, y)^{\nu-\gamma}}
 \m(\d y)
\end{align*}
and hence $|f|\d\m\in K_{\nu,\beta}^{\,p}$.

When $\nu=\beta$, note that the function
$r^{\nu-\gamma} \bigl(\log{r^{-1}}\bigr)^p$ is
monotonically increasing for sufficiently small $r>0$ and
converges to $0$ as $r\downarrow 0$.
Then we have for such $r$,
\begin{align*}
 \int_{B_r(x)}
  \bigl(\log{d(x, y)^{-1}}\bigr)^p |f(y)|
 \m(\d y)
&=
 \int_{B_r(x)}
  \frac{\,|f(y)|\,}{d(x, y)^{\nu-\gamma}}
  d(x, y)^{\nu-\gamma} \bigl(\log{d(x, y)^{-1}}\bigr)^p
 \m(\d y)
\\
&\leq
 r^{\nu-\gamma}(\log r^{-1})^p
 \int_{B_r(x)}
  \frac{\,|f(y)|\,}{d(x, y)^{\nu-\gamma}}
 \m(\d y)
\end{align*}
and hence $|f|\d\m\in K_{\nu,\beta}^{\,p}$. 
\end{pf}

\begin{apf}{Theorem~\ref{Th3}}
The assertion of Theorem~\ref{Th3} is a special case of 
the assertion of Theorem~\ref{Th4} by setting $\nu=\alpha$. 
So there is no need to show the proof. 
\end{apf}

\begin{apf}{Corollary~\ref{cor:Th3}}
Suppose $\nu-p(\nu-\beta)>0$. 
It suffices to prove that 
$d(\cdot,o)^{-\gamma}\in L^q_{\rm\tiny unif}(E;\m)$ for any 
$q\in ]\nu/(\nu-p(\nu-\beta)),\nu/\gamma[$ with $\nu\geq\beta$, and for $q\in [1,\nu/\gamma[$ with $\nu<\beta$. By 
{\bf (A)}\ref{asmp:Bishop}, there exists $C>0$ such that  
 $\sup_{x\in E}\m(B_r(x))\leq V(r)\leq 
Cr^{\nu}$ for any $r\in]0,2]$. Then
\begin{align*}
\sup_{d(x,o)\geq 2}\int_{d(x,y)\leq1}d(y,o)^{-q\gamma}\m(\d y)
&\leq \int_{d(x,y)\leq1,d(y,o)\geq 1}d(y,o)^{-q\gamma}\m(\d y)
\\
&\leq\sup_{x\in E}\m(B_2(x))<+\infty.
\end{align*}
On the other hand, 
\begin{align*}
\sup_{d(x,o)< 2}\int_{d(x,y)\leq1}d(y,o)^{-q\gamma}\m(\d y)
&\leq \int_{d(y,o)<3}d(y,o)^{-q\gamma}\m(\d y)
\\
&\leq
\sum_{k=0}^{\infty}\int_{\frac{\,3\,}{2^{k+1}}\leq d(y,o)<\frac{\,3\,}{2^k}}d(y,o)^{-q\gamma}\m(\d y)\\
&\leq 
2C\cdot3^{\nu-q\gamma}\sum_{k=0}^{\infty}\left(\frac{\,1\,}{2^{\nu-q\gamma}} \right)^{k}<+\infty.
\end{align*}
 Therefore
\begin{align*}
 \sup_{x\in E}
 \int_{d(x,y)\leq1}d(y,o)^{-q\gamma}\m(\d y)
<
 +\infty.
\end{align*}
\end{apf}

\subsection{Proofs of Theorem \ref{Th5}, Corollary \ref{cor:Th5} and Theorem \ref{Th6}}

\begin{apf}{Theorem~\ref{Th5}}
First we prove \ref{item:Th5:1}.
{\setlength{\leftmargini}{0mm}
\begin{enumerate}
\item[](Case I) $\nu<\beta$:
In this case, we can directly check
$\mu\in K_{\nu,\beta}^{\,p}$ by
$\sup_{x\in E}\mu(B_r(x))\leq C r^{\eta}$ for $r\in]0,r_0]$.

\item[](Case II) $\nu=\beta$: 
In this case, 
\begin{align}
\int_0^{\infty}u^{\nu-(1-\delta)\beta-1}
\Phi_2(u)\d u<+\infty
\label{eq:estimateIINT}
\end{align}
for any $\delta\in]0,1]$.
By the upper estimate, we have
\begin{align*}
 \int_0^tp_s(x,y)\d s
&\leq
 \frac{\,\beta\,}{d(x,y)^{\nu-\beta}}
 \int_{d(x,y)/t^{1/\beta}}^{\infty}
  u^{\nu-\beta-1}\Phi_2(u)
 \d u
\\
&\leq
 \frac{\,t^{\delta}\,}{d(x,y)^{\nu-(1-\delta)\beta}}
 \left(
  \beta
  \int_{0}^{\infty}
   u^{\nu-(1-\delta)\beta-1}\Phi_2(u)
  \d u
 \right)
\end{align*} 
and hence, for $r\in]0,r_0]$ and $t\in]0,t_0[$,
\begin{align}
\hspace{-1cm}\sup_{x\in E}&
\int_{d(x,y)\leq r}
 \left(\int_0^tp_s(x,y)\d s \right)^p
\mu(\d y)
\nonumber
\\
&\leq 
 t^{p\delta}
 \left(
  \beta
  \int_{0}^{\infty}
   u^{\nu-(1-\delta)\beta-1}\Phi_2(u)
  \d u
 \right)^p
 \sup_{x\in E}
 \int_{d(x,y)\leq r}
  \frac{\,1\,}{d(x,y)^{p(\nu-(1-\delta)\beta)}}
 \mu(\d y)
\nonumber
\\
&=
 t^{p\delta}
 \left(
  \beta
  \int_{0}^{\infty}
   u^{\nu-(1-\delta)\beta-1}\Phi_2(u)
  \d u
 \right)^p
 \sup_{x\in E}
 \sum_{k=0}^{\infty}
 \int_{\frac{\,r\,}{2^{k+1}}< d(x,y)\leq\frac{\,r\,}{2^k}}
  \frac{\,1\,}{d(x,y)^{p(\nu-(1-\delta)\beta)}}
 \mu(\d y)
\nonumber
\\
&\leq
 t^{p\delta}
 \left(
  \beta
  \int_{0}^{\infty}
   u^{\nu-(1-\delta)\beta-1}\Phi_2(u)
  \d u
 \right)^p
 C_2
 2^{p(\nu-\beta)+\delta p\beta}
 \left(
  \sum_{k=0}^{\infty}
  2^{-k(\eta-p(\nu-\beta)-\delta p\beta)}
 \right)
 r^{\eta-p(\nu-\beta)-\delta p\beta}.
\label{eq:INest}
\end{align}
This goes to $0$ as $r\to0$ for a   
constant $\delta\in]0,1[$ satisfying 
$\eta-p(\nu-\beta)>\beta p\delta$. Hence
$\mu\in K_{\nu,\beta}^{\,p}$ by Theorem~\ref{Th1}. 

\item[](Case III) $\nu>\beta$:
In this case, 
\eqref{eq:estimateIINT} holds 
for any $\delta\in[0,1]$. So \eqref{eq:INest}  goes to 
$0$ as $r\to0$ for a constant $\delta\in]0,1[$ satisfying 
$\eta-p(\nu-\beta)>\beta p\delta$. Hence 
$\mu\in K_{\nu,\beta}^{\,p}$ by Theorem~\ref{Th1}. 
\end{enumerate}
}

Next we prove \ref{item:Th5:2}.
Assume $\eta-p(\nu-\beta)<0$.
Then it implies $\nu>\beta$. 
Suppose that 
there exist $x_0\in E$ and $r_0,C_1>0$ such that 
$\mu(B_r(x_0))\geq C_1r^{\eta}$ for all $r\in]0,r_0]$.
Then, 
for any $s\in]0,+\infty[$ and $t\in]0,t_0]$, 
\begin{align}
 \sup_{x\in E}&
 \int_{d(x,y)\leq s}
  \left(\int_0^tp_s(x,y)\d s \right)^p
 \mu(\d y)
\nonumber
\\
&\geq
 \lim_{r\to0}
 \left(
  \beta
  \int_{r_0/t^{1/\beta}}^{\infty}
   u^{\nu-\beta-1}\Phi_1(u)
  \d u
 \right)^p
 \int_{d(x_0,y)\leq r}
  \frac{\,1\,}{d(x_0,y)^{p(\nu-\beta)}}
 \mu(\d y)
\nonumber
\\
&\geq
 \lim_{r\to0}
 \left(
  \beta
  \int_{r_0/t^{1/\beta}}^{\infty}
   u^{\nu-\beta-1}\Phi_1(u)
  \d u
 \right)^p
 \frac{\;\mu(B_r(x_0))\;}{r^{p(\nu-\beta)}}
% \sum_{k=0}^{\infty}
% \int_{\frac{\,r\,}{2^{k+1}}< d({\color{blue}{x_0}},y)
%\leq\frac{\,r\,}{2^k}}
%  \frac{\,1\,}{d({\color{blue}{x_0}},y)^{p(\nu-\beta)}}
% \mu(\d y)
\nonumber
\\
&\geq
 \lim_{r\to0}
 \left(
  \beta
  \int_{r_0/t^{1/\beta}}^{\infty}
   u^{\nu-\beta-1}\Phi_1(u)
  \d u
 \right)^p
 %\biggl(
 C_1
 % \sum_{k=0}^{\infty}
 %  2^{-k(\eta-p(\nu-\beta))}
 %\biggr)
 r^{\eta-p(\nu-\beta)}
=
 +\infty
.
\label{eq:INest*}
\end{align}
This implies $\mu\notin D_{\nu,\beta}^{\,p}$
under $\eta-p(\nu-\beta)<0$
by Theorem~\ref{Th1}
${\rm (1')}\Longleftrightarrow{\rm (5')}$.
Note here that
$
 \int_0^{\infty}u^{\nu-1}\Phi_1(u)\d u
\leq
 \int_0^{\infty}u^{\nu-1}\Phi_2(u)\d u
<
 +\infty
$.

Finally we prove \ref{item:Th5:3}.
Assume $\eta-p(\nu-\beta)\leq0$.
Then it implies $\nu>\beta$.
Suppose that
there exist $x_0\in E$ and $r_0,C_1,C_2>0$ such that 
$C_1r^{\eta}\leq \mu(B_r(x_0))\leq C_2r^{\eta}$
for all $r\in]0,r_0]$.
Take a large $\ell>0$ so that $C_1>C_2/\ell^{\eta}$.
Then, 
for any 
$s\in]0,+\infty[$,
$r\in]0,s\land r_0]$ and
$t\in]0,t_0]$, 
\begin{align}
 \sup_{x\in E}&
 \int_{d(x,y)\leq s}
  \left(\int_0^tp_s(x,y)\d s \right)^p
 \mu(\d y)
\nonumber
\\
&\geq
 \left(
  \beta
  \int_{r_0/t^{1/\beta}}^{\infty}
   u^{\nu-\beta-1}\Phi_1(u)
  \d u
 \right)^p
 \int_{d(x_0,y)\leq r}
  \frac{\,1\,}{d(x_0,y)^{p(\nu-\beta)}}
 \mu(\d y)
\nonumber
\\
&=
 \left(
  \beta
  \int_{r_0/t^{1/\beta}}^{\infty}
   u^{\nu-\beta-1}\Phi_1(u)
  \d u
 \right)^p
 \sum_{k=0}^{\infty}
 \int_{\frac{\,r\,}{\ell^{k+1}}< d(x_0,y)\leq\frac{\,r\,}{\ell^k}}
  \frac{\,1\,}{d(x_0,y)^{p(\nu-\beta)}}
 \mu(\d y)
\nonumber
\\
&\geq
 \left(
  \beta
  \int_{r_0/t^{1/\beta}}^{\infty}
   u^{\nu-\beta-1}\Phi_1(u)
  \d u
 \right)^p
 \left(
  C_1-\frac{C_2}{\ell^{\eta}}
 \right)
 \left(
  \sum_{k=0}^{\infty}
  \ell^{-k(\eta-p(\nu-\beta))}
 \right)
 r^{\eta-p(\nu-\beta)}
\nonumber
\\
&=
 +\infty.
\label{eq:INest**}
\end{align}
This implies
$\nu\notin D_{\nu,\beta}^{\,p}$ under $\eta-p(\nu-\beta)\leq0$
by Theorem~\ref{Th1}${\rm (1')}\Longleftrightarrow{\rm (5')}$.
\end{apf}

\begin{apf}{Corollary~\ref{cor:Th5}}
It suffices to consider $\eta=\nu$.
In view of 
\cite[Corollary~4.1]{KwT:Katounderheat},
if further {\bf (A)}\ref{asmp:SC} is satisfied,
then the Ahlfors regularity holds.
This implies the assertion. 
\end{apf}

\begin{apf}{Theorem~\ref{Th6}}
We will estimate
the inner integral
$
 \int_{d(x, y)\leq r}
 \left(
  \int_0^t
   p_s(x, y)
  \d s
 \right)^{\,p}
 \mu(\d y)
$
and the outer integral
$
 \int_{d(x, y)> r}
 \left(
  \int_0^t
   p_s(x, y)
  \d s
 \right)^{\,p}
 \mu(\d y)
$
respectively.
{\setlength{\leftmargini}{6mm}
\begin{enumerate}
\item  
Suppose $\mu(E)<+\infty$. 
Let
$
 0
<
 \delta
<
 \frac{\,\eta - p(\nu - \beta)\,}{p\beta}
\hspace{1mm}
 (\leq 1)
$
and fix $r\in ]0, r_0[$.
We have for $t<1\wedge t_0 \wedge r^\beta$
\begin{align*}
 \int_0 ^t
  s^{-\nu/\beta}
  \Phi_2(rs^{-1/\beta})
 \d s
=
 \beta
 \int_{r/t^{1/\beta}} ^\infty
  u^{\nu-\beta-1}
  \Phi_2(u)
 \d u
\leq
 \left\{
  \beta
  \int_{1} ^\infty
   u^{\nu-1}
   \Phi_2(u)
  \d u
 \right\}
 r^{-\beta} t
\end{align*}
and then
\begin{align}
\notag
 \sup_{x\in E}
 \int_{d(x, y)> r}
 \left(
  \int_0^t
   p_s(x, y)
  \d s
 \right)^{\,p}\!\!
 \mu(\d y)
&\leq
 \mu(E)
 \left\{
  \beta
  \int_{1} ^\infty
   u^{\nu-1}
   \Phi_2(u)
  \d u
 \right\}^{\,p}\!\!
 r^{-p\beta}
 t^p
\\
\label{eq:estp=1outer}
&\leq
 \mu(E)
 \left\{
  \beta
  \int_{1} ^\infty
   u^{\nu-1}
   \Phi_2(u)
  \d u
 \right\}^{\,p}\!\!
 r^{-p\beta}
 t^{p\delta}
.
\end{align}
As for the outer integral,
take $\xi >0$ such that
\begin{equation*}
 \frac{\,\beta-\nu\,}{\beta}
<
 \xi
\quad
 \text{and}
\quad
 \delta
<
 \xi
<
 \frac{\,\eta-p(\nu-\beta)\,}{p\beta}
\end{equation*}
(the first assertion automatically holds when $\nu\geq\beta$).
We note that
\eqref{eq:INest} with replacing $\delta$ with $\xi$
still holds in the case $\nu\not=\beta$
because of $(\beta-\nu)/\beta<\xi$.
Hence we have for $t<1\wedge t_0$
\begin{align}
\label{eq:estp=1inner}
 \sup_{x\in E}
 \int_{d(x, y)\leq r}
 \left(
  \int_0^t
   p_s(x, y)
  \d s
 \right)^{\,p}
 \mu(\d y)
\leq
 C
 r^{\eta-p(\nu-\beta)-p\xi\beta}
 t^{p\xi}
\leq
 C
 r^{\eta-p(\nu-\beta)-p\xi\beta}
 t^{p\delta}.
\end{align}
Therefore we conclude that $\mu\in S_K^{\,p,\delta}$.

\item 
Suppose $p>1$ and $\mu\in S_D$. 
The inner integral estimate \eqref{eq:estp=1inner}
still holds in this case, and
it remains to calculate the outer integral.
If
$p\geq \frac{\eta+\beta}{\nu} (>1)$
and
$\delta< \frac{\eta-p(\nu-\beta)}{p\beta}$,
we have
$p\delta\leq p-1$.
If
$p>1$
and
$\delta< \frac{p-1}{p}\cdot\frac{\eta-p(\nu-\beta)}{p\beta}$,
we also have
$p\delta \leq p-1$.
Hence, for both cases,
by a similar calculation
as that to obtain \eqref{eq:estp=1outer},
we have for $t<1\wedge t_0 \wedge r^\beta$
\begin{align}
\notag
 \sup_{x\in E}&
 \int_{d(x, y)> r}
 \left(
  \int_0^t
   p_s(x, y)
  \d s
 \right)^p
 \mu(\d y)
\\
\notag
&\leq
 \left\{
  \sup_{x\in E}
  \int_E
  \left(
   \int_0^1
    p_s(x, y)
   \d s
  \right)
  \mu(\d y)
 \right\}
 \left\{
  \beta
  \int_{1} ^\infty
   u^{\nu-1}
   \Phi_2(u)
  \d u
 \right\}^{\,p-1}\!\!
 r^{-(p-1)\beta}
 t^{p-1}
\\
\label{eq:estp=1outer_2}
&\leq
 \left\{
  \sup_{x\in E}
  \int_E
  \left(
   \int_0^1
    p_s(x, y)
   \d s
  \right)
  \mu(\d y)
 \right\}
 \left\{
  \beta
  \int_{1} ^\infty
   u^{\nu-1}
   \Phi_2(u)
  \d u
 \right\}^{\,p-1}\!\!
 r^{-(p-1)\beta}
 t^{p\delta}
.
\end{align}
Therefore we conclude that $\mu\in S_K^{\,p,\delta}$.

\item
Suppose $p>1$ and $\mu\in S_D$ and assume that 
there exists $C>0$ such that
$\sup_{x\in E}\mu(B_r(x))\leq Cr^{\eta}$
for all $r\in]0,+\infty[$. 
Under this condition, 
the calculation to obtain \eqref{eq:INest}
remains valid for all $r\in]0,+\infty[$. 
Let
$
 0
<
 \delta
<
 \frac{\,\eta - p(\nu - \beta)\,}{p\beta}
\hspace{1mm}
 (< 1)
$
and take $\xi >0$ such that
\begin{equation*}
 \frac{\,\beta-\nu\,}{\beta}
<
 \xi
,
\quad
 \delta
<
 \xi
<
 \frac{\,\eta-p(\nu-\beta)\,}{p\beta}
\quad
 \text{and}
\quad
 \frac{\,p\,}{p-1}\delta
<
 \frac
 {\eta-p(\nu-\beta)-p\beta\delta}
 {\eta-p(\nu-\beta)-p\beta\xi}
\end{equation*}
(the first assertion automatically holds when $\nu\geq\beta$).
Set $r=t^{-\alpha}$, where
$
 \alpha
=
 \frac{p(\xi-\delta)}{\eta-p(\nu-\beta)-p\beta\xi}
>
 0
$.
Then,
by the same calculation
as that to obtain \eqref{eq:estp=1outer_2},
we have for $t<1\wedge t_0$
\begin{align*}
 \sup_{x\in E}&
 \int_{d(x, y)> r}
 \left(
  \int_0^t
   p_s(x, y)
  \d s
 \right)^p
 \mu(\d y)
\\
&\leq
 \left\{
  \sup_{x\in E}
  \int_E
  \left(
   \int_0^1
    p_s(x, y)
   \d s
  \right)
  \mu(\d y)
 \right\}
 \left\{
  \beta
  \int_1 ^\infty
   u^{\nu-1}
   \Phi_2(u)
  \d u
 \right\}^{p-1}
 r^{(p-1)\beta}
 t^{p-1}
\\
&\leq
 \left\{
  \sup_{x\in E}
  \int_E
  \left(
   \int_0^1
    p_s(x, y)
   \d s
  \right)
  \mu(\d y)
 \right\}
 \left\{
  \beta
  \int_1 ^\infty
   u^{\nu-1}
   \Phi_2(u)
  \d u
 \right\}^{p-1}
 t^{p\delta},
\end{align*}
where we used
\begin{align*}
 (p-1)(\alpha\beta+1)
=
 (p-1) 
 \frac
 {\eta-p(\nu-\beta)-p\beta\delta}
 {\eta-p(\nu-\beta)-p\beta\xi}
>
 p\delta
.
\end{align*}
The inner integral estimate follows from
the middle side of \eqref{eq:estp=1inner}
since
\begin{equation*}
 r^{\eta-p(\nu-\beta)-p\xi\beta}
 t^{p\xi}
=
 t^{-\alpha(\eta-p(\nu-\beta)-p\xi\beta)}
 t^{p\xi}
=
 t^{p\delta}.
\end{equation*}
Therefore we conclude that $\mu\in S_K^{\,p,\delta}$.

We next assume that
$
 \int_1^\infty
  u^{\nu+\gamma-1}\Phi_2(u)
 \d u
<
 +\infty
$
holds for any $\gamma>0$.
The inner integral estimate \eqref{eq:estp=1inner}
still holds in this case.
As for the outer integral,
take $\gamma>0$ such that
$
 \frac{\beta+\gamma}{\beta}
>
 \frac{\!p\!}{p-1}\delta
$.
Then, by a similar calculation
as that to obtain \eqref{eq:estp=1outer},
we have
\begin{align*}
 \sup_{x\in E}&
 \int_{d(x, y)> r}
 \left(
  \int_0^t
   p_s(x, y)
  \d s
 \right)^p
 \mu(\d y)
\\
&\leq
 \left\{
  \sup_{x\in E}
  \int_E
  \left(
   \int_0^1
    p_s(x, y)
   \d s
  \right)
  \mu(\d y)
 \right\}
 \left\{
  \beta
  \int_{1} ^\infty
   u^{\nu+\gamma-1}
   \Phi_2(u)
  \d u
 \right\}^{\,p-1}\!\!
 r^{-(p-1)(\beta+\gamma)}
 t^{(p-1)\frac{\beta+\gamma}{\beta}}
\\
&\leq
 \left\{
  \sup_{x\in E}
  \int_E
  \left(
   \int_0^1
    p_s(x, y)
   \d s
  \right)
  \mu(\d y)
 \right\}
 \left\{
  \beta
  \int_{1} ^\infty
   u^{\nu-1}
   \Phi_2(u)
  \d u
 \right\}^{\,p-1}\!\!
 r^{-(p-1)(\beta+\gamma)}
 t^{p\delta}
.
\end{align*}
Therefore we conclude that $\mu\in S_K^{\,p,\delta}$.
\end{enumerate}
}
\vspace{-2eM}
\end{apf}

\section{Examples}

\begin{ex}[{Brownian motions on \boldmath$\R^d$}]\label{ex:Brownian}
{\rm
Let
${\bf X}^{\rm w}=(\Omega, B_t, {\PP}_x)_{x\in\R^d}$ be
a $d$-dimensional Brownian motion on $\R^d$. 
Consider $p\in[1,+\infty[$. 
We say that $\mu\in K_{d}^{\,p}$ (or $\mu\in K_{d,2}^{\,p}$) if and only if 
\begin{align*}
 \lim_{r\to0}
 \sup_{x\in\R^d}
 \int_{|x-y|<r}\frac{\,\mu(\d y)\,}{\;|x-y|^{(d-2)p}\;}
=0
\quad&
 \text{for }\quad d\geq3
,
\\
 \lim_{r\to0}
 \sup_{x\in\R^d}
 \int_{|x-y|<r}(\log|x-y|^{-1})^p\mu(\d y)
=0
\quad&
 \text{for }\quad d=2
,
\\
 \sup_{x\in\R^d}
 \int_{|x-y|\leq1}
 \mu(\d y)
<
 +\infty
\quad&
 \text{for }\quad d=1
.
\end{align*}
We write $K_d$ instead of $K_d^1$ for $p=1$. 
As in Section \ref{sec:introduction}, 
we have $K_d^{\,p}=S_K^{\,p}$ by 
\cite[Example~2.4]{TM:pKato} or its extension 
Theorem~\ref{Th2}. 

The $d$-dimensional Lebesgue measure $\m$ belongs to  
$K_d^{\,p}=S_K^{\,p}$
if and only if $p\in [1, d/(d-2)_+[$ by
Theorem~\ref{Th5} or Corollary~\ref{cor:Th3}, where
$d/(d-2)_+:=d/(d-2)$ if $d\geq3$,
$d/(d-2)_+:=+\infty$ if $d=1,2$. 
For any non-negative bounded $g\in L^1(\R^d)$
the finite measure $g\m$ also belongs to $S_K^{\,p,\delta}$
for $0<\delta<(d-p(d-2))/2p$ under $p\in[1,d/(d-2)_+[$
by Theorem~\ref{Th6}(1).
Moreover,
$\m\in S_K^{\,p,\delta}$
for $0<\delta<(d-p(d-2))/2p$ under $p\in]1,d/(d-2)_+[$
by Theorem~\ref{Th6}(3).

The surface measure $\sigma_R$
on the $R$-sphere $\partial B_R(0)$ satisfies that 
$\sigma_R(B_r(x))\leq C_2r^{d-1}$
for any $x\in\R^d$ and $r>0$ with some $C_2>0$, and 
$\sigma_R(B_r(x))\geq C_1r^{d-1}$
for any $x\in\partial B_R(0)$ and $r\in]0,r_0[$
with some $C_1,r_0>0$.
Then we can conclude that  
$\sigma_R\in K_{d}^{\,p}=S_K^{\,p}$ holds if and only if
 $p\in [1, (d-1)/(d-2)_+[$ under $d\geq2$
by Theorems~\ref{Th2} and \ref{Th5}, where
$(d-1)/(d-2)_+:=(d-1)/(d-2)$ if $d\geq3$,
$(d-1)/(d-2)_+:=+\infty$ if $d=2$. 
Moreover,
$\sigma_R\in S_K^{\,p,\delta}$ holds
for $0<\delta<((d-1)-p(d-2))/2p$ 
under $p\in]1, (d-1)/(d-2)_+[$ with $d\geq2$
by Theorem~\ref{Th6}(1).
%
%In particular,
%any Hausdorff measure with its Hausdorff dimension 
%{\color{blue}{
%$\eta\in]\,(d-\alpha)_+,\,d\,]$ (resp.~$\eta\in]\,(d-2)_+,\,d\,]$)}}
%is of $p$-Kato class in the framework of
%rotationally symmetric $\alpha$-stable processes
%(resp.~Brownian motions on Riemmanian manifolds),
%provided {\color{blue}{$p\in[\,1,\,\eta/(d-\alpha)_+[$
%(resp.~$\eta\in[\,1,\,\eta/(d-2)_+[$)}}.
%
By Theorem~\ref{Th3}, we also have that  
$f\in L^q_{\text{\tiny unif }}(\R^d)$ implies
$|f|\d\m\in S_K^{\,p}$ if $q>d/(d-p(d-2))$ with $d\geq2$,
or $q\geq1$ with $d=1$.
}
\end{ex}

\begin{ex}[{Symmetric relativistic \boldmath$\alpha$-stable process, symmetric \boldmath$\alpha$-stable process}]\label{ex:relative}
{\rm
{Take $0 <\alpha <2$ and $m\geq0$. 
Let $\text{\bf X}=(\Omega,X_t,{\PP}_x)$ be a  
L\'evy process on $\R^d$ with 
\begin{equation*}
{\EE}_{0}\left[e^{\sqrt{-1}\<  \xi, X_{t}  \>}\right]=
\exp \left(-t\left\{(|\xi|^2+m^{2/\alpha})^{\alpha/2}-m\right\}\right).
\end{equation*}
If $m>0$, it is called the \emph{relativistic $\alpha$-stable process with mass $m$} 
(see \cite{CKS:POTA2012}).  
In particular, if $\alpha=1$ and $m>0$,
it is called the \emph{free relativistic Hamiltonian process}
(see \cite{CMS:1990, CKS:AOP2012, HS:PP}). 
When $m=0$, $\text{\bf X}
$ is nothing but the usual \emph{(rotationally) symmetric 
$\alpha$-stable process}. 
It is known that $\text{\bf X}
$ is transient if and only if $d>2$ under $m>0$ or $d>\alpha$ under $m=0$, and 
{\bf X} is a doubly Feller conservative process.    

Let $(\E,\F)$ 
be the Dirichlet form on $L^2(\R^d)$ associated with $\text{\bf X}$. 
Using Fourier transform $\hat{f}(x):=\frac{\,1\,}{(2\pi)^{d/2}}\int_{\R^d}e^{i\<x,y\>}f(y){\rm d}y$, it follows 
from \cite[Example~1.4.1]{FOT} that 
\begin{align*}
\left\{\begin{array}{rl}
 \F
\hspace{-3mm}
&=
 \displaystyle{\left\{
  f\in L^2(\R^d)
 \;\left|\;
  \int_{\R^d}
   |\hat{f}(\xi)|^2
   \left((|\xi|^2+m^{2/\alpha})^{\alpha/2} -m \right)
  \d \xi
 <
  +\infty\right.
 \right\}}
,
 \\
 \E(f,g)
\hspace{-3mm}
&=
 \displaystyle{
 \int_{\R^d}
  \hat{f}(\xi)\bar{\hat{g}}(\xi)
  \left((|\xi|^2+m^{2/\alpha})^{\alpha/2} -m \right)
 \d \xi
\quad
 \text{ for }\quad f,g\in \F.
}
\end{array}\right.
\end{align*}
It is shown in \cite{CS3} that the corresponding jumping measure 
$J$ of $(\E,\F)$  
satisfies
\begin{align*}
J({\rm d}x{\rm d}y)=J_m(x,y){\rm d}x{\rm d}y \quad \text{ with }~~J_m(x,y)=A(d,-\alpha)
\frac{\,\Psi(m^{1/\alpha}|x-y|)\,}{|x-y|^{d+\alpha}},
\end{align*}
where $A(d,-\alpha)=\frac{\,\alpha 2^{d+\alpha}\Gamma(\frac{d+\alpha}{2})\,}{2^{d+1}\pi^{d/2}\Gamma(1-\frac{\,\alpha\,}{2})}$, and $\Psi(r):=I(r)/I(0)$ with 
\begin{equation*}
I(r):=\int_0^{\infty}s^{\frac{\,d+\alpha\,}{2}-1}e^{-\frac{\,s\,}{4}-\frac{\,r^2\,}{s}}{\rm d}s
\end{equation*}
is a decreasing function satisfying 
$\Psi(r)\asymp e^{-r}(1+r^{(d+\alpha-1)/2})$ near $r=+\infty$, and $\Psi(r)=1+\Psi''(0)r^2/2+o(r^4)$ near $r=0$. 
In particular, 
\begin{align*}
\left\{\begin{array}{rl}
 \F
\hspace{-3mm}
&=
 \displaystyle{\left\{
  f\in L^2(\R^d)
 \;\left|\;
  \int_{\R^d\times\R^d}
   |f(x)-f(y)|^2J_m(x,y)
  \d x \d y
 <
  +\infty\right.
 \right\}}
,
\\
 \E(f,g)
\hspace{-3mm}
&=
 \displaystyle{
 \frac{\,1\,}{2}
 \int_{\R^d\times\R^d}
  (f(x)-f(y))(g(x)-g(y))J_m(x,y)
 \d x \d y
\quad\text{ for }\quad f,g\in \F
}.
\end{array}\right.
\end{align*}
Let $p_t(x,y)$ be the heat kernel of {\bf X}.
The following global heat kernel estimate is proved in \cite[Theorem~2.1]{CKS:JFA2012}: 
there exists $C_1,C_2>0$ such that 
\begin{align}
C_2^{-1}\Phi^m_{1/C_1}(t,x,y)\leq p_t(x,y)\leq C_2\Phi^m_{C_1}(t,x,y),\label{eq:GHKERS}
\end{align}
where 
\begin{align*}
 \Phi^m_C(t,x,y)
:=
 \left\{
 \begin{array}{ll}
  t^{-d/\alpha} \wedge  tJ_m(x, y),
 & t\in]0,1/m], \\ 
  m^{d/\alpha-d/2} t^{-d/2}
  \exp\left(
   -C^{-1}(m^{1/\alpha}|x - y| \wedge m^{2/\alpha-1}\frac{\,|x-y|^2\,}{t})
  \right),
 & t\in]1/m,+\infty[.
 \end{array}
 \right.
\end{align*}
In particular, we have 
\begin{align}
C_2^{-1}(t^{
-d/\alpha} \wedge  tJ_m(x, y))\leq p_t(x,y)\leq C_2(t^{
-d/\alpha} \wedge  tJ_m(x, y))\quad\text{ for }\quad t\in]0,1/m]
.\label{eq:GHKES}
\end{align}
It is shown
in \cite[Theorem~1.2 and Example~2.4]{CK:MixedHeat} 
or \cite[Theorem~1.2]{CKK:TranAMS2011,CKK:TranAMS2015}
that $p_t(x,y)$ is jointly continuous in 
$(t,x,y)\in ]0,+\infty[\times \R^d\times \R^d$.
The $\beta$-order resolvent kernel
$r_{\beta}(x,y)\in[0,+\infty]$ is also continuous
in $(x,y)\in \R^d\times\R^d$.
}
Consider $p\in[1,+\infty[$. 
We say that $\mu\in K_{d,\alpha}^{\,p}$ if and only if 
\begin{align*}
 \lim_{r\to0}
 \sup_{x\in \R^d}
 \int_{|x-y|<r}\frac{\,\mu(\d y)\,}{\;|x-y|^{p(d-\alpha)}\;}
=
 0
\quad&
 \text{for }\quad d>\alpha
,
\\
 \lim_{r\to0}
 \sup_{x\in \R^d}
 \int_{|x-y|<r}\hspace{0cm}(\log |x-y|^{-1})^p\mu(\d y)
=
 0
\quad&
 \text{for }\quad d=\alpha \,(=1)
,
\\
 \sup_{x\in\R^d}
 \int_{|x-y|\leq1}\mu(\d y)
<
 +\infty
\quad&
 \text{for }\quad \alpha>d \,(=1)
.
\end{align*}
Then we have $K_{d,\alpha}^{\,p}=S_K^{\,p}$
by Theorem~\ref{Th2}. 

Consequently,
the $d$-dimensional Lebesgue measure $\m$ belongs to  
$K_{d,\alpha}^{\,p}=S_K^{\,p}$ if and only if
$\alpha>\frac{\,(p-1)d\,}{p}$
by Corollary~\ref{cor:Th3} or Theorem~\ref{Th5}, and
for any non-negative bounded $g\in L^1(\R^d)$
the finite measure $g\m$ also belongs to $S_K^{\,p,\delta}$
for $0<\delta<1-\frac{\,(p-1)d\,}{p\alpha}$
under $p\in]1,d/(d-\alpha)_+[$ 
by Theorem~\ref{Th6}(1).
Moreover, 
$\m\in S_K^{\,p,\delta}$
for $0<\delta<1-\frac{\,(p-1)d\,}{p\alpha}$
under $p\in]1,d/(d-\alpha)_+[$ by Theorem~\ref{Th6}(3). 
Here $d/(d-\alpha)_+:=d/(d-\alpha)$ if $d>\alpha$ and 
$d/(d-\alpha)_+:=+\infty$ if $d\leq\alpha$.  
The surface measure $\sigma_R$
on the $R$-sphere $\partial B_R(0)$ satisfies that 
$\sigma_R(B_r(x))\leq C_2r^{d-1}$
for any $x\in\R^d$ and $r>0$ with some $C_2>0$, and 
$\sigma_R(B_r(x))\geq C_1r^{d-1}$
for any $x\in\partial B_R(0)$ and $r\in]0,r_0[$
with some $C_1,r_0>0$.
Then we can conclude that  
$\sigma_R\in K_{d}^{\,p}=S_K^{\,p}$ holds if and only if
$\alpha>\frac{\,(p-1)d+1\,}{p}$ under $d>\alpha$
by Theorems~\ref{Th2} and \ref{Th5}, and
$\sigma_R\in S_K^{\,p,\delta}$ holds for
$0<\delta<1-\frac{\,(p-1)d+1\,}{p\alpha}$
under $p\in]1, (d-1)/(d-\alpha)[$ with $d\geq2$
by Theorem~\ref{Th6}(1). 
%Here $(d-1)/(d-\alpha)_+:=(d-1)/(d-\alpha)$ if $d>\alpha$, $(d-1)/(d-\alpha)_+:=+\infty$ if $d\leq\alpha$. 
%Since  
%the surface measure $\sigma_R$ on the $R$-sphere $\partial B_R(0)$ satisfies $C_1r^{d-1}\leq \sigma_R(B_r(x))\leq C_2r^{d-1}$ for any $x\in\R^d$ {\color{blue}{and $r\in]0,+\infty[$}} with some $C_1,C_2>0$, 
%we can conclude that  
%$\sigma_R\in K_{d,\alpha}^{\,p}=S_K^{\,p}$ holds if and only if  $\alpha>\frac{\,(p-1)d+1\,}{p}$ {\color{blue}{under $d\geq2$}}
%by Theorems~\ref{Th2} and \ref{Th5}, 
%and $\sigma_R\in  S_K^{\,p,\delta}$ holds for $0<\delta<1-\frac{\,(p-1)d+1\,}{p\alpha}$ {\color{blue}{
%under $p\in]1,(d-1)/(d-\alpha)_+[$ with $d\geq2$ 
%by Theorem~\ref{Th6}(1).}} 
%Moreover, since any Hausdorff measure $\mathcal{H}^{\eta}$ on a subset of $\R^d$ with its Hausdorff dimension $\eta$
%satisfies $C_1r^{\eta}\leq \mathcal{H}^{\eta}(B_r(x))\leq C_2r^{\eta}$ ({\color{blue}{$r\in]0,+\infty[$}}) for some $C_1,C_2>0$, we can conclude $\mathcal{H}^{\eta}\in K_{d,\alpha}^{\,p}=S_K^{\,p}{\color{blue}{\subset K_{d,\alpha}=S_K}}$ holds if and only if $\alpha>\frac{\,pd-\eta\,}{p}$ 
%{\color{blue}{under $\eta>0$}}
%by Theorems~\ref{Th2} and \ref{Th5}. {\color{blue}{By Theorem~\ref{Th6}(1),(2), 
%$\mathcal{H}^{\eta}\in S_K^{\,p,\delta}$ holds for $0<\delta<1-\frac{\,pd-\eta\,}{2p}$ {\color{blue}{under 
%$p\in]1,\eta/(d-\alpha)_+[$ with $\eta>0$}} 
%provided 
%$\mathcal{H}^{\eta}(\R^d)<+\infty$ or $p>1$ under $\alpha p>pd-\eta$. 
%
By Theorem~\ref{Th3}, we also have that   
$f\in L^q_{\text{\tiny unif }}(\R^d)$ implies
$|f|\d\m\in S_K^{\,p}$ if
$q>d/(d-p(d-\alpha))$ with $d\geq\alpha$, or
$q\geq1$ with $d<\alpha$. 

We finally expose propositions on $K_{d,\alpha}^{\,p}$: 
\begin{prop}\label{prop:rearrangement}
Let $f$ be a $[0,+\infty]$-valued function on $[0,+\infty[$. 
Suppose that $|V(x)|\m(\d x)\in K_{d,\alpha}^{\,p}$ with $V(x):=f(|x|)$.  
Then we have 
\begin{align*}
 \int_0^R
  r^{d-p(d-\alpha)-1}f(r)
 \d r
<
 +\infty
&\quad 
 \text{ for \ some }\quad R>0\quad \text{ if }\quad d>\alpha
,
\\ 
 \int_0^R
  (\log(r^{-1}))^pf(r)
 \d r
<
 +\infty
&\quad 
 \text{ for \ some }\quad R\in]0,1/e\,[\quad \text{ if }\quad d=\alpha=1
,
\\   
 \int_0^Rf(r)\d r
<
 +\infty
&\quad 
 \text{ for \ some }\quad R>0\quad \text{ if }\quad \alpha>d=1.
\end{align*}
If further $f$ is decreasing on $[0,+\infty[$ and vanishes at infinity,  
then the converse holds.  
\end{prop}

\begin{pf}
Suppose $|V(x)|\m(\d x)\in K_{d,\alpha}^{\,p}$. From Theorem~\ref{Th2}, we have 
\begin{equation*}
 \sup_{x\in\R^d}\int_{|x-y|<R}|V(y)|\d y
<
 +\infty
\end{equation*}
for any $R\in]0,+\infty[$. Then we see that for any $R\in]0,+\infty[$ ($R\in]0,1/e[$ if $d=\alpha=1$),  
\begin{align*}
 \sup_{x\in\R^d}
 \int_{|x-y|<R}\frac{\,|V(y)|\,}{\;|x-y|^{p(d-\alpha)}\;}\d y
<
 +\infty
\quad&
 \text{if }\quad d>\alpha,
\\ 
 \sup_{x\in\R^d}
 \int_{|x-y|<R}(\log|x-y|^{-1})^p|V(y)|\d y
<
 +\infty
\quad&
 \text{if }\quad d=\alpha(=1),
\\ 
 \sup_{x\in\R^d}
 \int_{|x-y|<R}|V(y)|\d y
<
 +\infty
\quad&
 \text{if }\quad \alpha>d=1.
\end{align*} 
Hence we have the assertion.
Suppose the converse with the decrease of $f$. 
Then the symmetric decreasing rearrangement $V^*$ of $V$ equals to $|V|$ (see Chapter 3 in Lieb-Loss \cite{LiebLoss:analysis}). The simplest rearrangement inequality (see \cite[Theorem~3.4]{LiebLoss:analysis}) 
tells us that 
\begin{align*}
\sup_{x\in \R^d}\int_{|x-y|<r}\frac{\,|V(y)|\,}{\;|x-y|^{p(d-\alpha)}\;}\d y&=\int_{|y|<r}\frac{\,|V(y)|\,}{\;|y|^{p(d-\alpha)}\;}\d y
\\&=(d\cdot\omega_d)\int_0^rs^{d-p(d-\alpha)-1}f(s)\d s \quad \text{ if }\quad d>\alpha, \\
\sup_{x\in \R^d}\int_{|x-y|<r}\log(|x-y|^{-1})^p|V(y)|\d y&=\int_{|y|<r}(\log|y|^{-1})^p|V(y)|\d y\\&=
2\int_0^r(\log s^{-1})^pf(s)\d s \quad \text{ if }\quad d=\alpha=1,  
\end{align*}
which tends to $0$ as $r\to0$, respectively. 
Here $\omega_d$ is the volume of the unit ball $B_1(0)$. 
We also have 
\begin{align*}
 \sup_{x\in \R^d}
 \int_{|x-y|<R}|V(y)|\d y=
 \int_{|y|<R}|V(y)|\d y=
 2\int_0^Rf(s)\d s
<
 +\infty
\quad \text{ if }\quad \alpha>d=1. 
\end{align*}
Then $|V(x)|\m(\d x)\in K_{d,\alpha}^{\,p}$. 
\end{pf}

\begin{prop}\label{prop:superset}
 Let $V$ be a measurable function satisfying that 
\begin{align*}
 \int_a^{\infty}
  \m(|V|\geq t)^{\frac{\,d-p(d-\alpha)\,}{d}}
 \d t
&<
 +\infty
\quad\text{ for \ some }\quad a>0
\quad\text{ if }\quad d>\alpha, \\ 
 \int_a^{\infty}
  F(\m(|V|\geq t))
 \d t
&<
 +\infty
\quad\text{ for \ some }\quad a>0
\quad\text{ if }\quad d=\alpha=1, \\   
 \int_a^{\infty}
  \m(|V|\geq t)
 \d t
&<
 +\infty
\quad\text{ for \ some }\quad a>0
\quad \text{ if }\quad \alpha>d=1.
\end{align*}
Then $|V(x)|\m(\d x)\in K_{d,\alpha}^{\,p}$. Here $\m$ is the Lebesgue measure on $\R^d$ and $F$ is a function defined by 
$F(s):=\int_0^{s/2}(\log^+ u^{-1})^p\d u$ with 
$\log^+u:=(\log u)\lor0$. 
\end{prop}

\begin{pf}
The proof is similar with \cite[Theorem~4.12]{AizSim:Brown}. 
We only prove the case $d>\alpha$. 
We may assume $|V|=V^*$, where $V^*$ is the symmetric decreasing rearrangement of $V$, 
because $\m(|V|>t)=\m(V^*>t)$ for all $t>0$. Hence, there exists a decreasing function 
$f$ on $]0,+\infty[$ such that $V(x)=f(|x|)$ and $V$ is lower semi-continuous.  
Let $f^{-1}(t):=\sup\{s>0\mid f(s)>t\}$ be the 
right continuous inverse of $f$, which is also a decreasing function. 
We may assume $0<f^{-1}(t)<+\infty$ for any $t>0$, that is, 
$f$ has an infinite limit at origin and 
no positive limit at infinity, because $a$ specified in the condition 
can be taken to be arbitrarily large. When $f$ has a finite limit at origin, $V$ is essentially 
bounded, which implies $|V(x)|\m(\d x)\in K_{d,\alpha}^{\,p}$. 
First we assume that $f$ is continuous, but we do not assume the strict decrease of $f$.
The continuity of $f$ yields that $b=f(f^{-1}(b))$ for any $b>0$ and we see 
 $f^{-1}(f(r))\leq r$ for any $r>0$. 
Then for $a<A$
\begin{align*}
\int_{f^{-1}(a)}^{f^{-1}(A)}r^{d-p(d-\alpha)}\d f(r)&\leq \int_{f^{-1}(a)}^{f^{-1}(A)}(f^{-1}(f(r)))^{d-p(d-\alpha)}\d f(r) \\ 
&=\int_{f(f^{-1}(a))}^{f(f^{-1}(A))}(f^{-1}(t))^{d-p(d-\alpha)}\d t\\
&=\int_a^A(f^{-1}(t))^{d-p(d-\alpha)}\d t\\ 
&\leq \int_a^{\infty}(f^{-1}(t))^{d-p(d-\alpha)}\d t<+\infty,
\end{align*}
because $\m(|V|\geq t)=\omega_d(f^{-1}(t))^d$. 
By way of the integration by parts formula for Riemann-Stieltjes integrals,   
\begin{align}
 (d-p(d-\alpha))&
 \int_0^{f^{-1}(a)}f(r)r^{d-p(d-\alpha)-1}\d r
\nonumber
\\
&\leq
 a(f^{-1}(a))^{d-p(d-\alpha)}
+
 \int_a^{\infty}(f^{-1}(t))^{d-p(d-\alpha)}\d t
<+\infty.
\label{eq:integrationpart}
\end{align}
Next we show (\ref{eq:integrationpart}) 
for general $f$. Note that $V^*(x)=f(|x|)$ is lower semi-continuous 
(lower semi-continuity of $f$ is clarified later). We set 
$f_n(t):=\inf\{f(s)+n|s-t|\mid s\in[0,+\infty[\}$. Then $\{f_n\}$ is an increasing sequence of 
 nonnegative $n$-Lipschitz function on $[0,+\infty[$.
We then see that 
\begin{equation*}
 f_n(|x|)
=
 \inf
 \left\{
  f(|z|)+n\bigl||z|-|x|\bigr|
 \hspace{1pt}\left|\hspace{3pt}
  z\in\R^d
 \right.
 \right\}.
\end{equation*}
Hence 
$f_n(|x|)$ converges to $V^*(x)=f(|x|)$ as $n\to\infty$, 
because of the lower semi-continuity of $V^*$,  
consequently, $f_n$ converges to $f$ as $n\to\infty$, hence $f$ is lower semi-continuous. 
Indeed, we may consider the case $f(|x|)>0$ and 
suppose the existence of $\varepsilon>0$ such that  $f_n(|x|)<f(|x|)-\varepsilon(>0)$ for all 
$n\in\N$. 
Then there exists $z_n\in\R^d$ with $f(|z_n|)+n||z_n|-|x||<f(|x|)-\varepsilon$. 
From this, we see $|z_n|\to|x|$ as $n\to\infty$ and may assume the existence of 
$z\in\R^d$ with $z_n\to z$ as $n\to\infty$ by taking a subsequence. 
Hence, we obtain a contradiction $f(|x|)=f(|z|)\leq\liminf_{n\to\infty}f(|z_n|)
\leq f(|x|)-\varepsilon$. 
We set $g_n(t):=\inf_{s\in[0,t]}f_n(s)$. Then $g_n$ is a decreasing continuous function 
vanishing at infinity. We see that 
$\uparrow\lim_{n\to\infty}g_n=f$. 
We also have that 
$\{g_n^{-1}\}$ is an increasing sequence and converges to $f^{-1}$ as $n\to\infty$ at each point. 
Since (\ref{eq:integrationpart}) holds
for $g_n,g_n^{-1}$, it holds for $f,f^{-1}$.
Therefore
the simplest rearrangement inequality shows 
that for $r>0$
\begin{align*}
 \sup_{x\in\R^d}
 \int_{|x-y|<r}\frac{\,|V(y)|\,}{\;|x-y|^{p(d-\alpha)}\;}\d y
\leq 
 \int_{|y|<r}\frac{\,|V^*(y)|\,}{\;|y|^{p(d-\alpha)}\;}\d y
=
 (d\cdot\omega_d)
 \int_0^rs^{d-p(d-\alpha)-1}f(s)\d s
\end{align*}
and the right-hand side goes to zero as $r\to 0$.
\end{pf}

\begin{cor}\label{cor:G}
 Let $d>\alpha$ and $G$  a positive increasing  
function on $]0,+\infty[$ satisfying 
\begin{equation*}
 \int_a^{\infty}
  (G'(s))^{1-\frac{\,d\,}{p(d-\alpha)}}
 \d s
<
 +\infty
\quad\text{ for \ some }\quad a>0.
\end{equation*}
Suppose that $\int_{\R^d}G(|V(x)|)\m(\d x)<+\infty$. Then $|V(x)|\m(\d x)\in K_{d,\alpha}^{\,p}$. 
\end{cor}

\begin{pf} The proof is quite similar as in \cite[Corollary~4.13]{AizSim:Brown}. 
We omit it. 
\end{pf}

\begin{rem}
{\rm 
The assertions in 
Propositions~\ref{prop:rearrangement}, \ref{prop:superset} and Corollary~\ref{cor:G} for $\alpha=2$ remain valid in the framework of $d$-dimensional Brownian motion 
${\bf X}^{\rm w}$
in Example~\ref{ex:Brownian}. 
}
\end{rem}
}
\end{ex} 

\begin{ex}[Jump type processes on \boldmath$d$-sets, cf. Chen-Kumagai~\cite{CK:HS}]
{\rm A metric measured space $(F,d,\m)$ satisfying the Ahlfors regularity that there exists $C>0$ such that 
$C^{-1}r^d\leq \m(B_r(x))\leq Cr^d$ for any $x\in F$ and $r\in]0,1[$ is called a \emph{$d$-set} if $F$ is a closed subset of 
$\R^n$ with $0<d\leq n$. For $\alpha\in]0,2[$, consider the following Dirichlet form $(\E^{(\alpha)},\F^{(\alpha)})$:
\begin{align*}
\left\{\begin{array}{rl}
 \F^{(\alpha)}
\hspace{-3mm}
&:=
 \left\{
  u\in L^2(F;\m)
 \;\left|\;
  \displaystyle{\int_F\int_F\frac{\,|u(x)-u(y)|^2\,}{|x-y|^{d+\alpha}}\m(\d x)\m(\d y)<+\infty}
 \right.
 \right\}
,
\\
 \E^{(\alpha)}(u,v)
\hspace{-3mm}
&:=
 \displaystyle{\frac12\int_F\int_F\frac{\,(u(x)-u(y))(v(x)-v(y))\,}{|x-y|^{d+\alpha}}\m(\d x)\m(\d y)}
,
 \quad 
 u,v\in\F^{(\alpha)}.
\end{array}\right.
\end{align*}
 
Under the condition that  for some $C>0$ with 
\begin{align}
\m(B_r(x))\leq Cr^d\quad\text{ for \ all }\quad x\in F\quad \text{ and }\quad r>0,\label{eq:Ahlforsupper}
\end{align}
Chen-Kumagai~\cite{CK:HS} shows that
the jump type process associated with 
$(\E^{(\alpha)},\F^{(\alpha)})$ admits a semigroup kernel
possessing the following upper and lower estimates: 
there exist $C_i=C_i(\alpha,d)>0$, $i=1,2$ such that
for all $(t,x,y)\in]0,1]\times F\times F$
\begin{equation*}
\hspace{-1cm}
 \frac{\,C_1\,}{t^{d/\alpha}}
 \frac1{
  \left(
   1+\frac{\,|x-y|\,}{t^{1/\alpha}}
  \right)^{d+\alpha}
 }
\leq
 p_t(x,y)
\leq
 \frac{\,C_2\,}{t^{d/\alpha}}
 \frac1{
  \left(
   1+\frac{\,|x-y|\,}{t^{1/\alpha}}
  \right)^{d+\alpha}
 }.
\end{equation*}
Then our conditions {\bf (A)}\ref{asmp:SC}, {\bf (A)}\ref{asmp:Bishop} and {\bf (A)}\ref{asmp:Phi} are satisfied in this context. 
Hu-Kumagai~\cite{HK:Nash} extends the result in \cite{CK:HS} replacing the embedding condition into $\R^n$ 
by a condition on the extension to a metric space with scaling property, however, they assume stronger assumption 
that for some $C>0$ with 
\begin{align}
C^{-1}r^d\leq \m(B_r(x))\leq Cr^d\quad\text{ for \ \ all }\quad x\in F\quad\text{ and }\quad r\in]0,\text{\rm diam}\,(F)].
\label{eq:Ahlforsupperlower}
\end{align}
One can apply our 
results under (\ref{eq:Ahlforsupperlower}). 
}
\end{ex}

\begin{ex}[Riemannian manifolds with lower Ricci curvature bound] 
{\rm \quad 
Let $(M,g)$ be a $d$-dimensional 
smooth complete Riemannian manifold with
${\Ric}_M\geq(d-1)\kappa_1$, $\kappa_1\in\R$. 
Let $\m$ be the volume measure induced from $g$ and
set $V(x,r):=\m(B_r(x))$. 
Under ${\Ric}_M\geq(d-1)\kappa_1$,
we have that the 
Bishop inequality $V(x,r)\leq V_{\kappa_1}(r)$ and
the Bishop-Gromov inequality 
$V(x,R)/V_{\kappa_1}(R)\leq V(x,r)/V_{\kappa_1}(r)$,
$0<r<R$ hold. 
Consequently,
we have the volume doubling condition
$\sup_{x\in M}V(x,2r)/V(x,r)<+\infty$ and  
$\int_1^{\infty}\frac{\,s\d s\,}{\log V(x,s)}=+\infty$
which implies    
the stochastic completeness of the 
Brownian motion on $(M,g)$ (see \cite{Gri:stoch}). 
Here
$V_{\kappa}(r)$ defined by 
$V_{\kappa}(r):=c_d\int_0^rS_{\kappa}(s)^{d-1}\d s$ 
with
$S_{\kappa}(s):=\frac{\,\sin t\sqrt{\kappa}\,}{\sqrt{\kappa}}$ if $\kappa>0$, $S_0(s)=s$, 
$S_{\kappa}(s):=\frac{\,\sinh t\sqrt{-\kappa}\,}{\sqrt{-\kappa}}$ if $\kappa<0$, 
is the volume of the ball with radius $r$ 
in the space form of constant sectional curvature $\kappa$ and 
$c_d$ is the volume of the unit ball in $\R^d$. 
So the condition {\bf (A)}\ref{asmp:SC} holds. 
We also have the scale invariant weak Poincar\'e inequality 
(depending on $\kappa_1$ if $\kappa_1<0$)
(see Saloff-Coste \cite{Salof:unif} or \cite[Theorem~5.6.5]{Salof:aspect}),
which implies the weak form of the weak Poincar\'e inequality
(see \cite[Theorem~5.5.1(i)]{Salof:aspect}). 
Then the heat kernel 
$p_t(x,y)$ of $(M,g)$ satisfies the following Li-Yau's estimate
(see \cite[Theorems~5.5.1 and 5.5.3]{Salof:aspect},  
cf. \cite[Theorems~6.1 and 6.2]{Gri:weightedmani}):
for each $T>0$ there exist $C_i=C_i(T)>0$, $i=1,2,3,4$ such that  
for $(t,x,y)\in ]0,T[\times M\times M$
\begin{equation*}
 \frac{\,C_3\,}{V(y,\sqrt{t})}
 \exp\left(-C_1\frac{\,d(x,y)^2\,}{t}\right)
\leq 
 p_t(x,y)
\leq
 \frac{\,C_4\,}{V(y,\sqrt{t})}
 \exp\left(-C_2\frac{\,d(x,y)^2\,}{t}\right).
\end{equation*}
The Bishop inequality tells us that 
{\bf (A)}\ref{asmp:Bishop} holds. 
Further
we assume that the injectivity radius of $M$
(write $\text{\rm inj}_M$) is positive. 
Then we have the following
(see the proof of \cite[Lemma~5]{Hebey:optimal} and
\cite[Proposition 14]{Croke:Iso}.
Though the framework of \cite{Croke:Iso} is restricted to
compact Riemannian manifolds,
the argument in \cite{Hebey:optimal} remains valid): 
there exists $C_d\in]0,+\infty[$ such that for any 
$r\in]0,\text{\rm inj}_M/2[$ and $x\in M$,  
\begin{equation*}
 V(x,r)\geq C_d r^d.
\end{equation*}
Hence for a small time $t_0>0$,
we have the Nash-type estimate:
for any $t\in]0,t_0[$,  
$
 \sup_{x,y\in E}p_t(x,y)
\leq
 C_4 t^{-d/2}
$,
which gives a Sobolev inequality under $d\geq3$ 
(see \cite{Hebey:optimal} again).  
Then \cite[Theorem~2.1]{KwT:Katofunc} holds. 
Also we have for any $t\in]0,t_0[$, $x,y\in M$ 
\begin{equation*}
 \frac{\,C_3\,}{t^{d/2}}
 \exp\left(-C_1\frac{\,d(x,y)^2\,}{t}\right)
\leq
 p_t(x,y)
\leq 
 \frac{\,C_4\,}{t^{d/2}}
 \exp\left(-C_2\frac{\,d(x,y)^2\,}{t}\right).
\end{equation*}
Since $]0,+\infty[\ni x\mapsto\sinh x/x$ is increasing,
$s\mapsto S_{\kappa_1}(s)/s^{d-1}$ is increasing
for $\kappa_1\leq0$, 
hence $s\mapsto V_{\kappa_1}(s)/s^d$ is so. 
For $\kappa_1>0$, $s\mapsto V_{\kappa_1}(s)/s^d$ is bounded.
We can confirm that 
\begin{equation*}
 \int_1^{\infty}
  \frac{\,(V_{\kappa_1}(s)\lor s^d)\,}{s}
  e^{-C_2 s^2}
 \d s
<
 +\infty.
\end{equation*}
Then {\bf (A)}\ref{asmp:Phi} holds. 
Therefore Theorems~\ref{Th2}, \ref{Th5}, \ref{Th1}, \ref{Th3}, \ref{Th4},  \ref{Th6} and Corollaries~\ref{cor:Th5}, \ref{cor:Th3} hold. 
In particular,
$\m\in S_K^{\,p}=K_{d}^{\,p}$ if and only if $d-p(d-2)>0$
by Corollary~\ref{cor:Th3}, and 
$\m\in S_K^{\,p,\delta}$ holds
for $\delta\in ]0,(d-p(d-2))/2p[$
under $p\in ]1,d/(d-2)_+[$
by the latter half of Theorem~\ref{Th6}(3).
}
\end{ex}

\begin{ex}[Nested fractals; cf.~\cite{Kuma:nested, FitzHamKuma:nested}]
{\rm The heat kernel of diffusion processes on the unbounded nested fractal $\tilde{K}$ 
constructed by 
Kumagai \cite{Kuma:nested} 
has the following upper and lower estimates: there exist $C_i>0$, $i=1,2,3,4$ such that 
for any $(t,x,y)\in]0,+\infty[\times \tilde{K}\times \tilde{K}$  
\begin{equation*}
\hspace{0cm}
 \frac{\,C_3\,}{t^{\frac{\,d_f\,}{d_w}}}
 \exp{\left(
  -C_1
  \left(
   \frac{\,d(x,y)\,}{t^{\frac{\,1\,}{d_w}}}
  \right)^{\frac{\,d_w\,}{d_J-1}}
 \right)}
\leq
 p_t(x,y)
\leq
 \frac{\,C_4\,}{t^{\frac{\,d_f\,}{d_w}}}
 \exp{\left(
  -C_2
  \left(
   \frac{\,d(x,y)\,}{t^{\frac{\,1\,}{d_w}}}
  \right)^{\frac{\,d_w\,}{d_J-1}}
 \right)}.
\end{equation*}
Here $d_f$ is the Hausdorff dimension of $\tilde{K}$, $d_w$ is called the 
\emph{walk dimension}, $d_J$ is a different constant from $d_w$. 
We consider $p\in[1,+\infty[$. 
It is known that $d_J=d_w$ if $\tilde{K}$ is the (unbounded) Sierpi\'nski Gasket. In general $d_f<d_w$.  
Hence $\mu\in S_K^{\,p}$ if and only if $\sup_{x\in \tilde{K}}\mu(\{y\in\tilde{K}:d(x,y)\leq1\})<+\infty$ and 
$L^q_{\text{\tiny unif }}(\tilde{K};\mu_{\tilde{K}})\subset S_K^{\,p}$ for $q\geq1$, where $\mu_{\tilde{K}}$ is the Hausdorff measure 
on $\tilde{K}$. In particular, $\mu_{\tilde{K}}\in S_K^{\,p}$ for 
$p\in[1,+\infty[$. 
}
\end{ex}

\begin{ex}[Sierpi\'nski Carpet; cf. \cite{BB:BS}]
{\rm On the unbounded Sierpi\'nski Carpet $\tilde{F}$, 
the heat kernel of diffusion processes exists and admits the following upper and lower 
estimates proved by Barlow-Bass \cite{BB:BS}: 
there exist constants $C_i>0$, $i=1,2,3,4$ such that 
for any $(t,x,y)\in]0,+\infty[\times \tilde{F}\times \tilde{F}$  
\begin{equation*}
 \frac{\,C_3\,}{t^{\frac{\,d_f\,}{d_w}}}
 \exp{\left(-C_1\left(\frac{\,|x-y|\,}{t^{\frac{\,1\,}{d_w}}}\right)^{\frac{\,d_w\,}{d_w-1}}\right)}
\leq
 p_t(x,y)
\leq
 \frac{\,C_4\,}{t^{\frac{\,d_f\,}{d_w}}}
 \exp{\left(-C_2\left(\frac{\,|x-y|\,}{t^{\frac{\,1\,}{d_w}}}\right)^{\frac{\,d_w\,}{d_w-1}}\right)}. 
\end{equation*}
Here $|x-y|$ denotes the Euclidean norm of $x-y$ in $\R^n$,  
$d_f$ is the Hausdorff dimension of $\tilde{F}$, $d_w$ is called the \emph{walk dimension} 
and $d_s:=2d_f/d_w$ is called the \emph{spectral dimension} of $\tilde{F}$. 
They have the relation  $1<d_s\leq d_f<n$, where $n$ is the dimension of the Euclidean space in 
which $\tilde{F}$ is embedded. 
Thus we have $2\leq d_w\leq 2n$. 

Take $p\in[1,+\infty[$. 
We say that $\mu\in K_{d_f,d_w}^{\,p}$ if and only if 
\begin{align*}
 \lim_{r\to0}
 \sup_{x\in \tilde{F}}
 \int_{\{y\in \tilde{F}: |x-y|<r\}}
 \frac{\,\mu(\d y)\,}{\;|x-y|^{p(d_f-d_w)}\;}
=
 0
,
\quad& \quad d_s>2
,
\\
 \lim_{r\to0}
 \sup_{x\in \tilde{F}}
 \int_{\{y\in \tilde{F}:|x-y|<r\}}
  \hspace{0cm}(\log |x-y|^{-1})^p
 \mu(\d y)
=
 0
,
\quad&
 \quad d_s=2
,
\\
 \sup_{x\in  \tilde{F}}
 \int_{\{y\in \tilde{F}:|x-y|\leq1\}}
 \mu(\d y)
<
 +\infty
,
\quad&\quad 
 d_s<2
.
\end{align*}
Let $\mu_{\tilde{F}}$ be the Hausdorff measure on $\tilde{F}$. 
In this case, Ahlfors regularity holds in the following sense that  
there exists $C>0$ such that $C^{-1}r^{d_f}\leq\mu_{\tilde{F}}(B_r(x))\leq Cr^{d_f}$ for all $r\in]1,+\infty[$ 
(see \cite[Lemma~2.3(f)]{BB:BS}).    
Then {\bf (A)}\ref{asmp:SC} is satisfied by 
\cite[Lemma~2.1]{KwT:Katounderheat}. 
By \cite[Remark~2.1]{KwT:Katounderheat},  
{\bf (A)}\ref{asmp:Bishop} holds by taking 
$V(r):=Cr^{d_f}$. We then have $K_{d_f,d_w}^{\,p}=S_K^{\,p}$ and 
$L^q_{\text{\tiny unif }}(\tilde{F};\mu_{\tilde{F}})\subset S_K^{\,p}$ 
if $q>d_s/(d_s-p(d_s-2))$ with $d_s\geq 2$, or $q\geq1$ with $d_s<2$. In particular, $\mu_{\tilde{F}}\in S_K^{\,p}$ for 
$p\in[1,+\infty[$ with $d_s-p(d_s-2)>0$. 
}
\end{ex}

\providecommand{\bysame}{\leavevmode\hbox to3em{\hrulefill}\thinspace}
\providecommand{\MR}{\relax\ifhmode\unskip\space\fi MR }
% \MRhref is called by the amsart/book/proc definition of \MR.
\providecommand{\MRhref}[2]{%
  \href{http://www.ams.org/mathscinet-getitem?mr=#1}{#2}
}
\providecommand{\href}[2]{#2}

\end{document}